\documentclass[12pt]{article}

\usepackage{url}
\usepackage{pstricks}
\usepackage{multirow}

\usepackage[font=sf,labelfont={sf,bf},scriptsize,margin=0.5in,labelsep=period]{caption}

\usepackage[parfill]{parskip}

\usepackage{graphicx}
\usepackage{amsmath}
\usepackage{amsthm}
\usepackage{float}
\usepackage{natbib}
\usepackage{geometry}
\usepackage{xspace}

\usepackage{hyperref}
\hypersetup{colorlinks=true,linkcolor=black,citecolor=black,urlcolor=black}

\newcommand{\titlecaption}[2][\relax]{\caption{\textbf{\sffamily{#1}}\newline \sffamily{#2}}}
 
\newcommand{\gsyn}{\ensuremath{g_{\text{syn}}}\xspace}
\newcommand{\Iapp}{\ensuremath{I_{\text{app}}}\xspace}

\newcommand{\Oh}[1]{\ensuremath{\mathcal{O}\!\left(#1\right)}\xspace}

\newcommand{\Ot}{\ensuremath{\Oh{t}}}

\newcommand{\given}[2]{\ensuremath{#1\left|#2\right.}}

\title{Conditioned Likelihoods Using Bifurcation Continuation in Inverse Modeling of Dynamical Systems\footnote{This research is supported in part by NSA grant H98230-11-1-0222 and NSF grant DMS-1062817.}}
\author{Karleigh Cameron\footnote{Central Michigan University, \protect\href{malito:camer1kj@cmich.edu}{\protect\nolinkurl{camer1kj@cmich.edu}}} \and Marissa Saladin\footnote{Aquinas College,\protect\href{malito:marissa.saladin@aquinas.edu}{\protect\nolinkurl{marissa.saladin@aquinas.edu}}}}

\begin{document}

\maketitle

\begin{abstract}

\noindent 
The Morris-Lecar (ML) model has applications to neuroscience and cognition. A simple network consisting of a pair of synaptically coupled ML neurons can exhibit a wide variety of deterministic behaviors including asymmetric amplitude state (AAS), equal amplitude state (EAS), and steady state (SS). In addition, in the presence of noise this network can exhibit mixed-mode oscillations (MMO), which represent the system being stochastically driven between these behaviors. In this paper, we develop a method to specifically estimate the parameters representing the coupling strength (\gsyn) and the applied current (\Iapp) of two reciprocally coupled and biologically similar neurons. This method employs conditioning the likelihood on cumulative power and mean voltage. Conditioning has the potential to improve the identifiability of the estimation problem. Conditioning likelihoods are typically much simpler to model than the explicit joint distribution, which several studies have shown to be difficult or impossible to determine analytically. We adopt a rejection sampling procedure over a closed defined region determined by bifurcation continuation analyses. This rejection sampling procedure is easily embedded within the proposal distribution of a Bayesian Markov chain Monte Carlo (MCMC) scheme and we evaluate its performance. This is the first report of a Bayesian parameter estimation for two reciprocally coupled Morris-Lecar neurons, and we find a proposal utilizing rejection sampling reduces parameter estimate bias relative to naive sampling.  Application to stochastically coupled ML neurons is a future goal.
\end{abstract}

\newpage
\tableofcontents

\newpage
\section{Introduction}

Transmembrane voltage is often recorded during physiological study of biological neurons.  However, voltage-gated ion channel activity and neurotransmitter levels are quite difficult to measure directly and are usually unobserved in such studies. In addition, there is a great diversity of neuron morphology, protein expression, and plasticity which may affect voltage dynamics and synaptic transmission~\citep{decarli,kollins}. Early development and senescence 
may also be major determinants of voltage response profiles~\citep{yeoman:insights,liu:gating}. Synaptic tuning in particular is thought to be an essential mediator of learning, stimulus response integration, and memory. There is evidence that memory and learning may depend critically on several distinct types of dynamic behavior in the voltage of neurons. 

The ML model reproduces the voltage of a single neuron and, depending on parameterization and initial conditions, can exhibit many of the experimentally observed behaviors of biological neurons~\citep{morrislecar}. In this paper, we explore a simple neural network consisting of two biologically identical, reciprocally coupled ML neurons. ~\citet{kuske} have shown that this modest model can exhibit a wide range of oscillating or non-oscillating voltage depending on the values of just a few parameters, specifically in this study, \Iapp and \gsyn. In the absence of noise, the model can predict synchronous or asynchronous firing, as well as either equal or unequal action potential amplitudes. Additionally, in the presence of even small noise in the applied current and weak synaptic coupling,  the system can exhibit mixed-mode oscillations (MMO) characterized by periods of small amplitude oscillation interrupted by large amplitude excursions.

In further work with the two ML neuron model, ~\citet{thompson:stochastic} explored two synaptically decoupled neurons driven by both common and independent intrinsic noise terms.  They found that shared common noise promotes synchronous firing of the two neurons, while separate intrinsic noise terms promoted asynchronous firing.  The relative scaling of the two noise sources was observed to be key in predicting the degree of synchrony.   In addition, while they did not specifically look at MMO, they hypothesized that such synchrony in a synaptically coupled network would increase the probability of MMO, by facilitating longer residence times within the unstable periodic orbits adjacent to the system's stable periodic orbits.  Indeed, in this paper we will detail the relative positions of these parameter regions as they are of key importance to our conditioned likelihood approach.  Specifically, we will provide a quick look-up table for the region in parameter space where  stable periodic orbits are possible.

\citet{Ditlevsen2012} develop a expectation-maximization (EM) stochastic particle filter method to estimate the parameters in a single ML neuron based on observation of voltage only.  A key aspect of their approach is that they assume both the voltage and the channel gating variables are in an oscillatory regime, but stochastically perturbed.  These perturbations are considered nuisance parameters which their method marginalizes away.  Specifically, they treat the unobserved channel gating variable from the model as a completely latent variable. Starting from estimates of the initial conditions for the voltage and channel gating variables, they iteratively predict the gating variable and voltage and then update the predicted voltage to the next time step using a modification of the well-known Euler differential equation solver. They discuss that an assumption of stationarity  in their method limits applicability to only short time windows over which current input can be considered constant (e.g. 600ms). They also note that certain parameters, conductances and reversal potentials in particular, are sensitive to choice of tuning parameters required by the method.

These studies demonstrate the active progress as well as the challenges of model parameter estimation for biological neuronal models and, more generally, for relaxation oscillator models.  Each of these studies derives asymptotic approximations or general forms for model likelihood, but use fundamentally different techniques and assumptions in doing so.  In each study the approach is specifically crafted to the model.  In this paper we attempt to develop a convenient Bayesian estimation scheme with only a few tuning parameters and relatively few mild assumptions. We focus our attention on deterministic synaptically coupled ML neurons.  Application of our method to stochastically coupled ML neurons is on-going  work in our group.

In the case of ML, estimation of \Iapp and \gsyn is non-trivial due to the diversity of possible dynamic behavior and the abrupt transitions among these seen with just small changes in these parameters' values.  However, we can better understand the critical values of these parameters by studying the system's bifurcation structure. We are able to locate parameter regimes where dramatic changes in the system appear. The neurons analyzed in this study are classified as Type II neurons, characterized by discontinuous drastic shifting between behavioral states. Because there is a distinct switch in behavior, bifurcation analyses determine a closed region of parameter space over which the relevant dynamics may occur. Sampling over such a feasibility region amounts to conditioning the inference on an \emph{a priori} assumed class of dynamics (e.g. stable node, limit cycle, steady state etc.).  While facilitating conditioning the likelihood on feature statistics of the voltage, this may translate into increased confidence and reduced bias in the parameter estimates. 

\section{Reciprocally Coupled ML Model}

\newcommand{\sinf}[1]{\ensuremath{\text{s}_{\infty}\left(#1\right)}}
\newcommand{\minf}[1]{\ensuremath{\text{m}_{\infty}\left(#1\right)}}
\newcommand{\lamb}[1]{\ensuremath{\lambda\left(#1\right)}}
\newcommand{\winf}[1]{\ensuremath{\text{w}_{\infty}\left(#1\right)}}

Our goal is parameter inference based on the temporal voltage response of two synaptically coupled neurons which are deterministically coupled to voltage-gated ionic conductance dynamics~\citep{morrislecar}. A single ML model has a two-dimensional phase space and is known to reproduce many of the behaviors experimentally observed in biological neurons~\citep{kuske}.  Therefore, systems of coupled ML neurons may offer a reasonable starting point for developing statistical inference methods for models of neuronal networks. The ML network we study is,
\begin{align*}
\frac{dv_1}{dt}& = \frac{1}{C} (g_{\text{Ca}} \cdot \minf{v_1} \cdot (v_1-v_{\text{Ca}})-g_\text{K} \cdot w_1(v_1-v_\text{K}) - g_\text{L} \cdot w_1(v_1-v_\text{L})\\ &+ \Iapp -\gsyn\cdot s_1\cdot (v_1 -v_{\text{syn}})) + \delta \cdot \xi_{1}\\
\frac{dv_2}{dt}& = \frac{1}{C} (g_{\text{Ca}} \cdot \minf{v_2} \cdot (v_2-v_{\text{Ca}})-g_\text{K} \cdot w_2(v_2-v_\text{K}) - g_\text{L} \cdot w_1(v_2-v_\text{L})\\ &+ \Iapp -\gsyn\cdot s_2\cdot (v_2 -v_{\text{syn}})) + \delta \cdot \xi_{2}\\
\frac{ds_{1}}{dt}& = \frac{\sinf{v_2} - s_{1}}{\tau}\\
\frac{ds_{2}}{dt}& = \frac{\sinf{v_1} - s_{2}}{\tau}\\
\frac{dw_{1}}{dt}& = \lamb{v1}(\winf{v1}-w1)\\
\frac{dw_{2}}{dt}& = \lamb{v2}(\winf{v2}-w2)
\end{align*}
where
\begin{align*}
	\minf{x}& = \frac{1}{2}\left(1+\tanh\left(\frac{x-v11}{v22}\right)\right)\\
	\sinf{x}& = \frac{1}{1+e^{-\frac{x-vt}{vs}}}\\
	\winf{x}& = \frac{1}{2}\left(1+\tanh\left(\frac{x-v3}{v4}\right)\right)\\
	\lamb{x}& = \phi\cosh{\frac{x-v3}{2v4}}.
\end{align*}

Note that in the stochastic version of this model $\delta>0$ and $\xi_{1}$ and $\xi_{2}$ are standard independent Wiener process variables.  In this paper, however we will be concerned with the deterministic version of this model where $\delta\equiv{0}$.

\begin{table}[H]\scriptsize
\begin{center}
\titlecaption[Parameters of ML model]{Assumed values of these parameters are given in the Appendix. Initial conditions of the variables may also determine the observed dynamics.}
\begin{tabular}{|c c || c c |}
\hline
\textbf{Variable} & \textbf{Definition} & \textbf{Variable} & \textbf{Definition} \\ \hline \hline
$I _{\text{app}}$ & Applied current  & $v_1, v_2, v_3, v_4$ & Membrane potentials \\
$g _{\text{syn}}$ & Coupling strength &  $v_{\text{L}}$, $v_{\text{Ca}}$, $v_\text{K}$, $v_{\text{syn}}$ & Equilibrium potentials\\
$C$ & Membrane capacitance & $\text{min}_f$ & $(1 + \tanh[(V-V_1)/V_2])/2$\\
$g_{\text{L}}$, $g_{\text{Ca}}$, $g_\text{K}$ & Conductance of membrane channels & $s_1, s_2$ & $\frac{1}{2} \cdot (1+ \tanh[\frac{v-v_3}{v_4}])$ \\
$w_1, w_2$ & Recovery variables &&\\
\hline
\end{tabular}
\label{tableofparams}
\end{center}
\end{table}

While the many parameters of the ML model are, in principle, experimentally verifiable, they impart a high dimension parameter space determining the behavior of the system. Indeed, there is a non-trivially large diversity of behaviors already possible by varying just a small few of these parameters. To simplify our exposition  and make the essence of our approach clear, the only two parameters we examine in detail are \Iapp, the exogenously applied current, and \gsyn, the synaptic coupling strength. Hence, in this study, \Iapp and \gsyn are the only unknown parameters and other parameters are assumed to have the values given in the Appendix. In order to estimate these parameters, we develop a Bayesian MCMC method based on a Metropolis-Hastings sampling approach with conditioned likelihood and rejection sampling proposal distribution.  

\section{Bifurcation Continuation Analyses}
Often a slight change in the parameters of a system causes a characteristic alteration in the system's behavior. The study of this change in behavior is known as bifurcation analysis. Information about the behavior of the ML model can be obtained by studying bifurcation diagrams.  In particular, parameter values corresponding to steady state and periodic solutions can be found. This information can then be used to better approximate \Iapp and \gsyn using MCMC.

\begin{figure}[H]
\begin{center}
\begin{tabular}{crlrl}
&a) & & b) & \\
\rotatebox{90}{\phantom{XXXXX}\large{Type I}}&& \includegraphics[scale=.3]{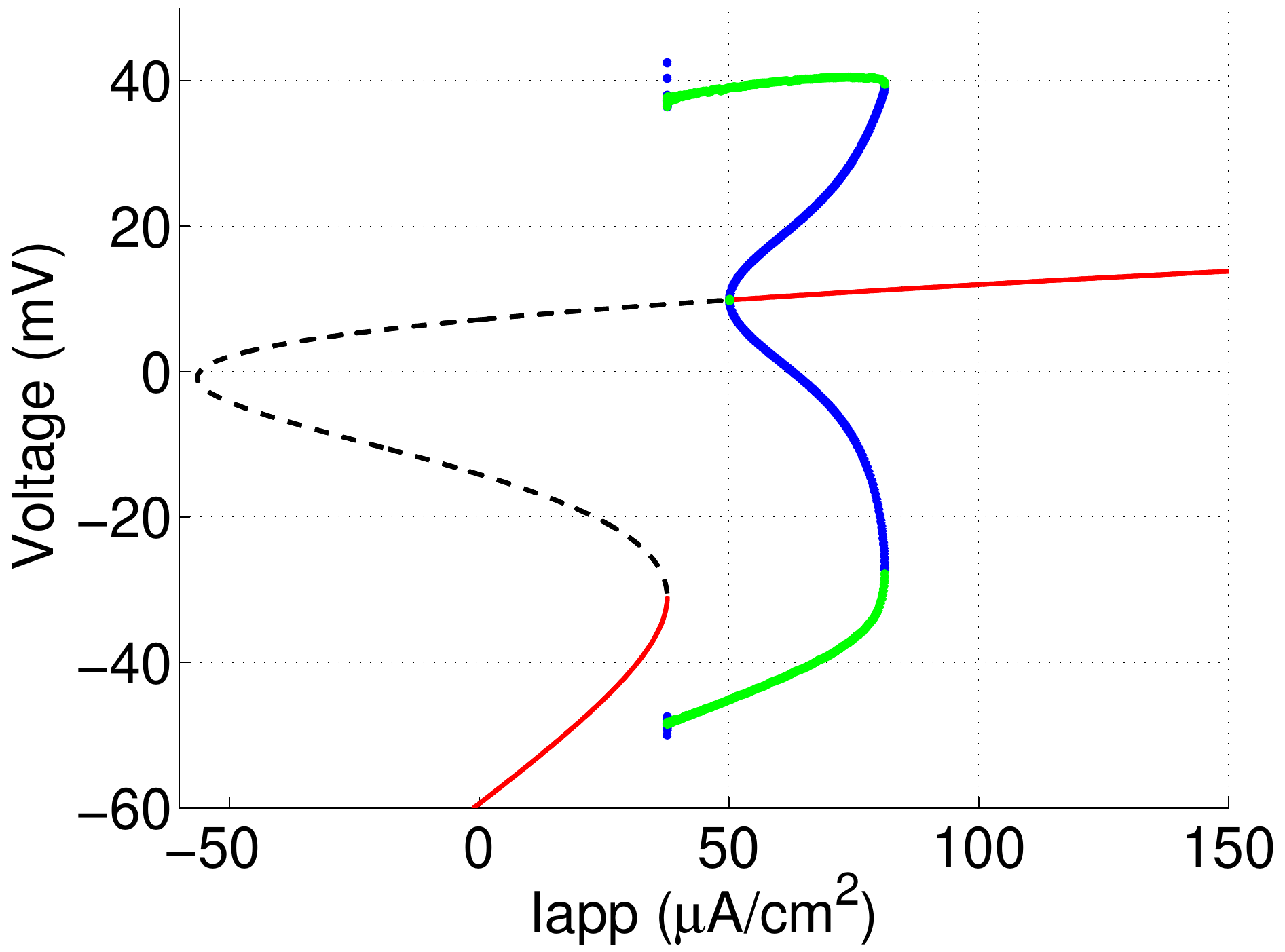} & & \includegraphics[scale=.3]{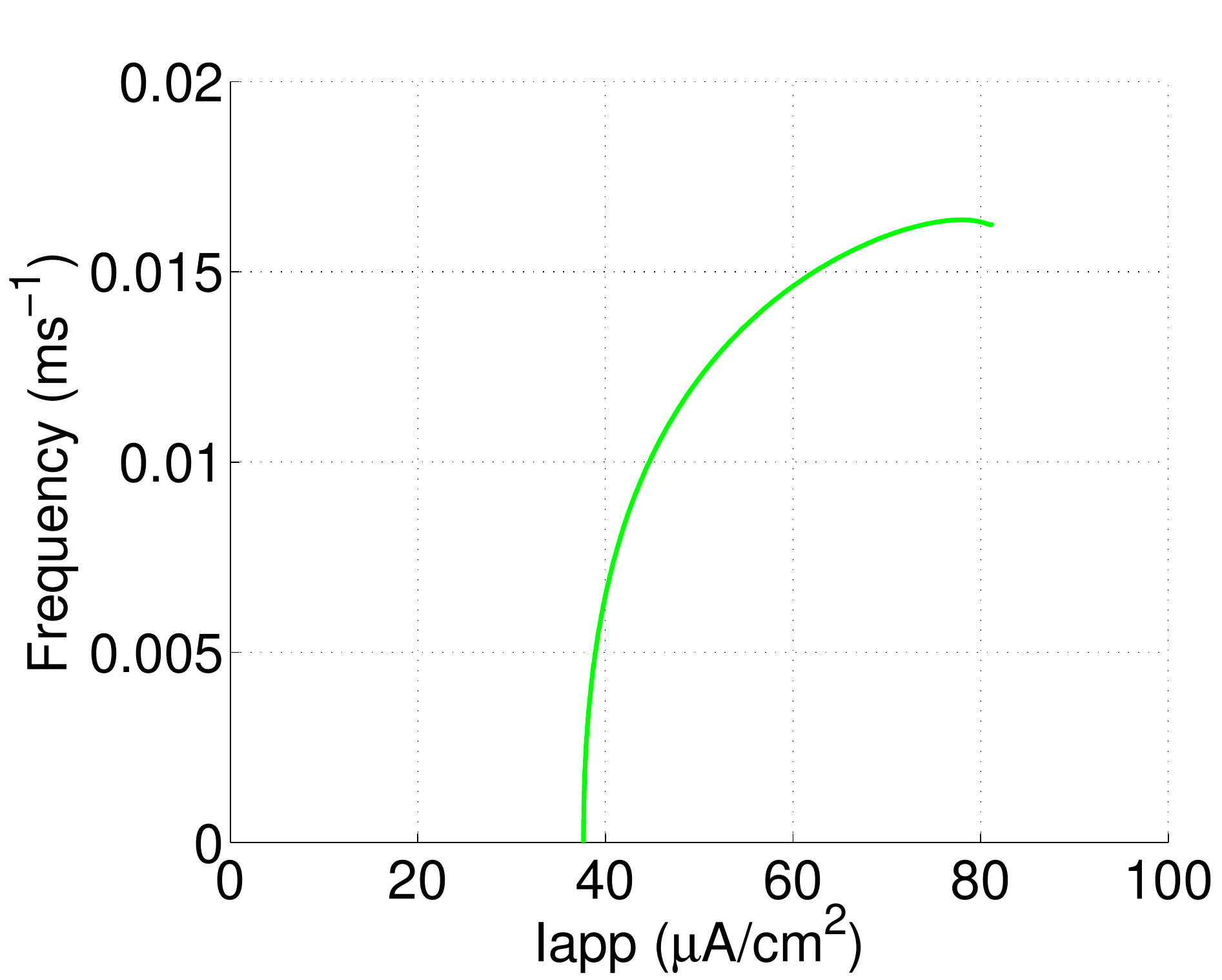}\\
&c) &  & d)\\
\rotatebox{90}{\phantom{XXXXX}\large{Type II}}&& \includegraphics[scale=.3]{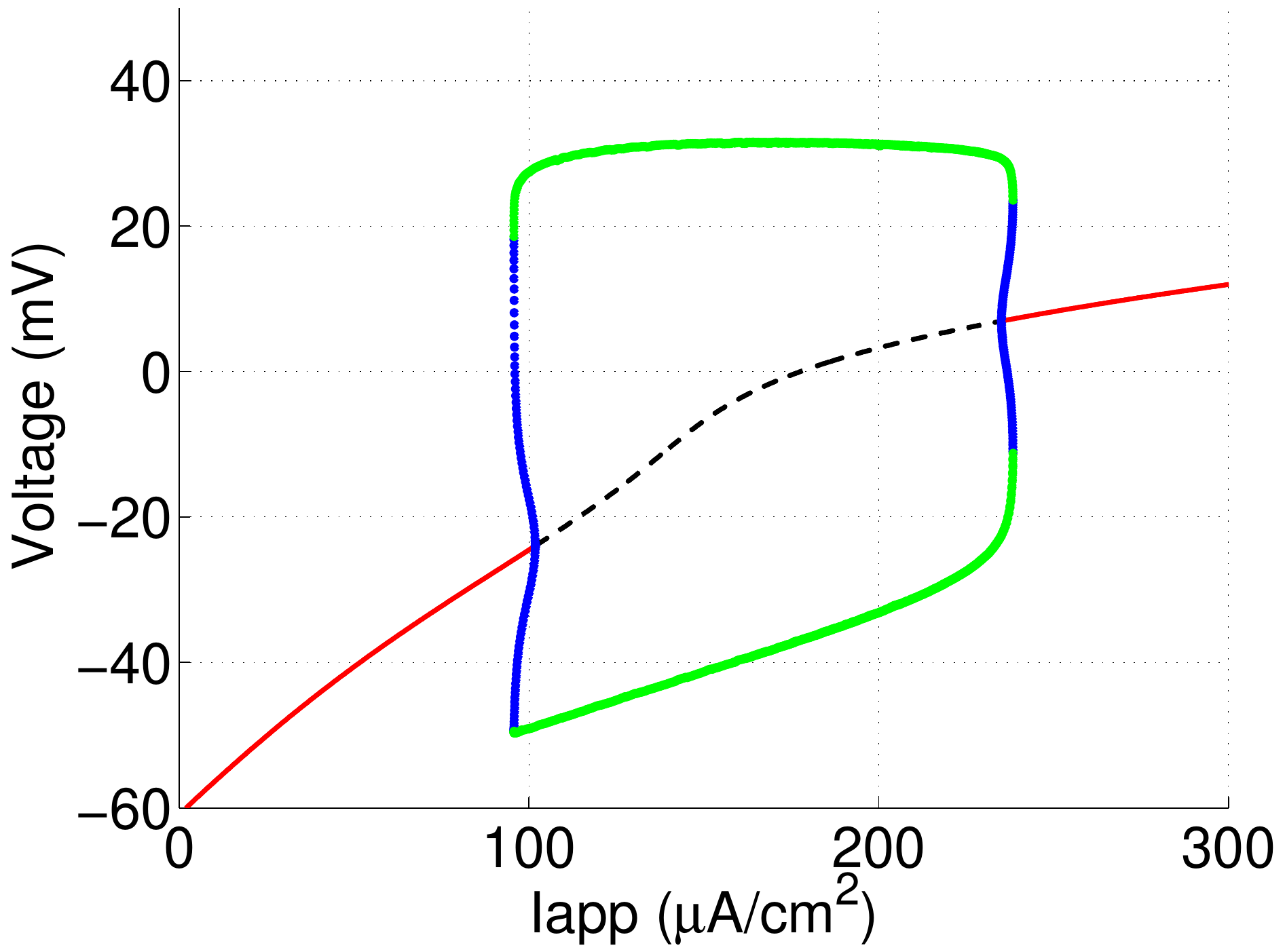} & & \includegraphics[scale=.3]{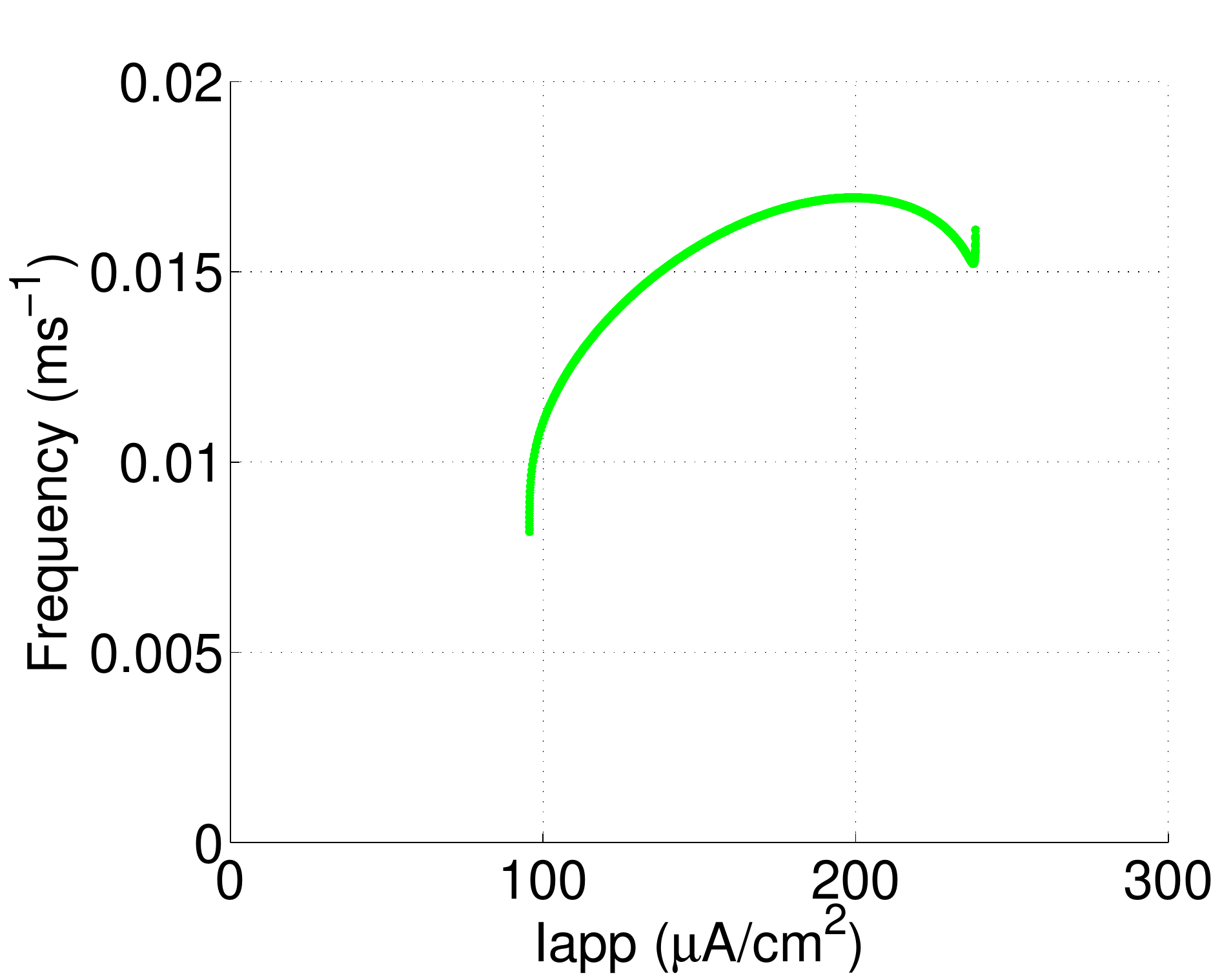}
\end{tabular}
\titlecaption[Type I vs Type II neurons]{Steady states are black, with a solid line depicting a stable state and a dashed line representing an unstable state. Oscillatory states are either green or blue. Green indicates a stable periodic state while blue represents an unstable periodic state.  All were generated with a \gsyn value of 0 to simulate uncoupled behavior. In a)--b), a Type I ML neuron is illustrated where all parameters are same as in Appendix except $v_{3}=15$ and $v_{4}=15$.   As seen in b), increasing \Iapp continuously increases the firing frequency from zero. In c)--d), a Type II ML neuron is illustrated where all parameters are identical to those in Appendix. In contrast to a Type I neuron, there is a discontinuous jump in frequency seen in d).}  
\label{N}
\end{center}
\end{figure}

For a single neuron, the bifurcation from steady state to oscillation is diagnostic for two broad classes referred to as {Type I} neurons and {Type II} neurons. In {Type I} neurons (Figs.~\ref{N}a) and \ref{N}b)), the switch from a {SS} to an oscillatory state is gradual and can often be difficult to detect precisely. However, in {Type II} neurons (Figs.~\ref{N}c) and \ref{N}d)), the change in state is sudden and drastic with a discontinuous jump in the frequency of action potentials. The neurons used in our study are determined to be Type II neurons. The sudden switch in behavior allows us to refine a proposal distribution for MCMC based on candidate parameter values.

\begin{figure}[H]
\begin{center}
\begin{tabular}{rlrl}
a) & & b) & \\
& \includegraphics[scale=.228]{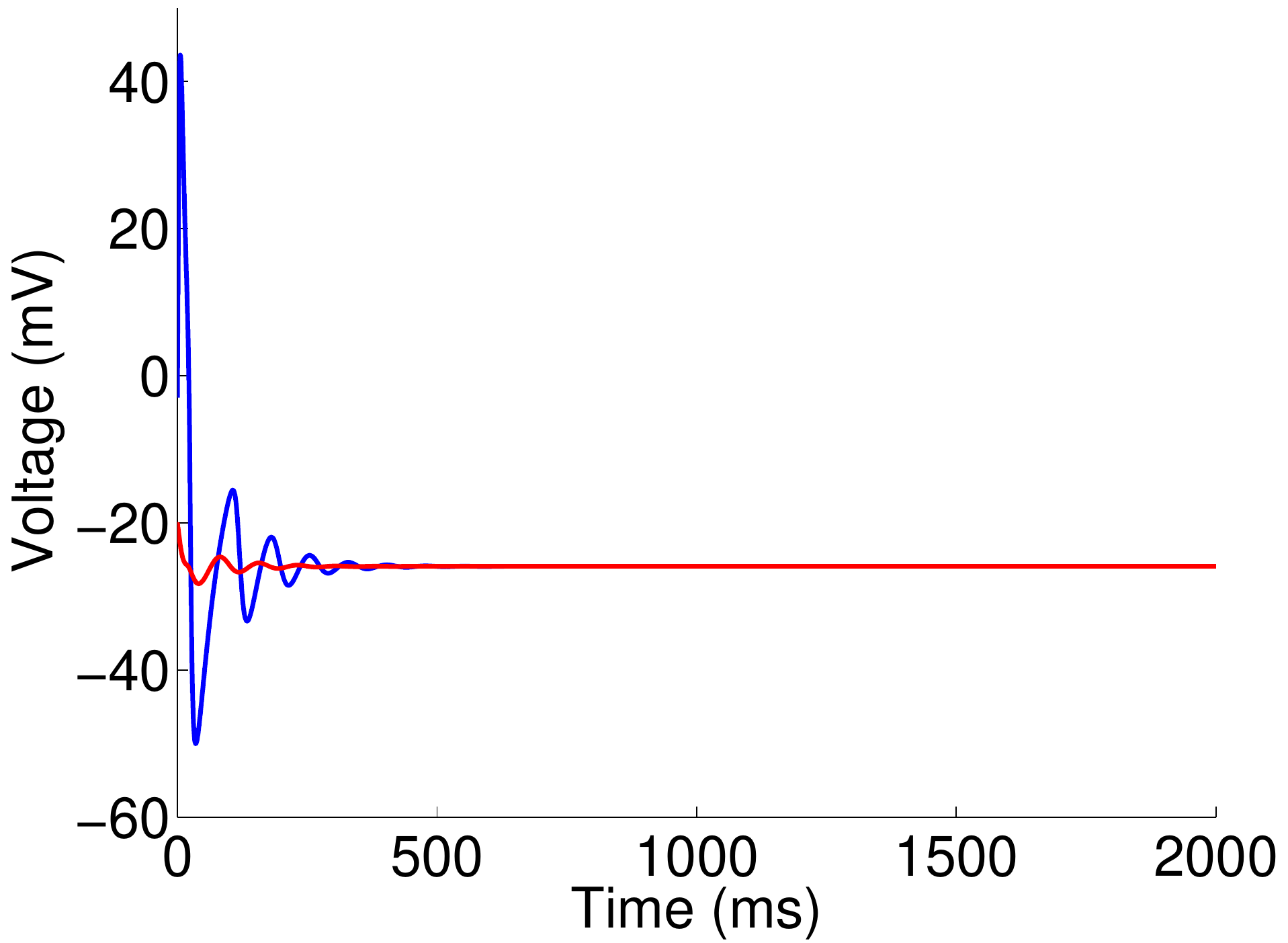}
\label{SS} & &\includegraphics[scale=.228]{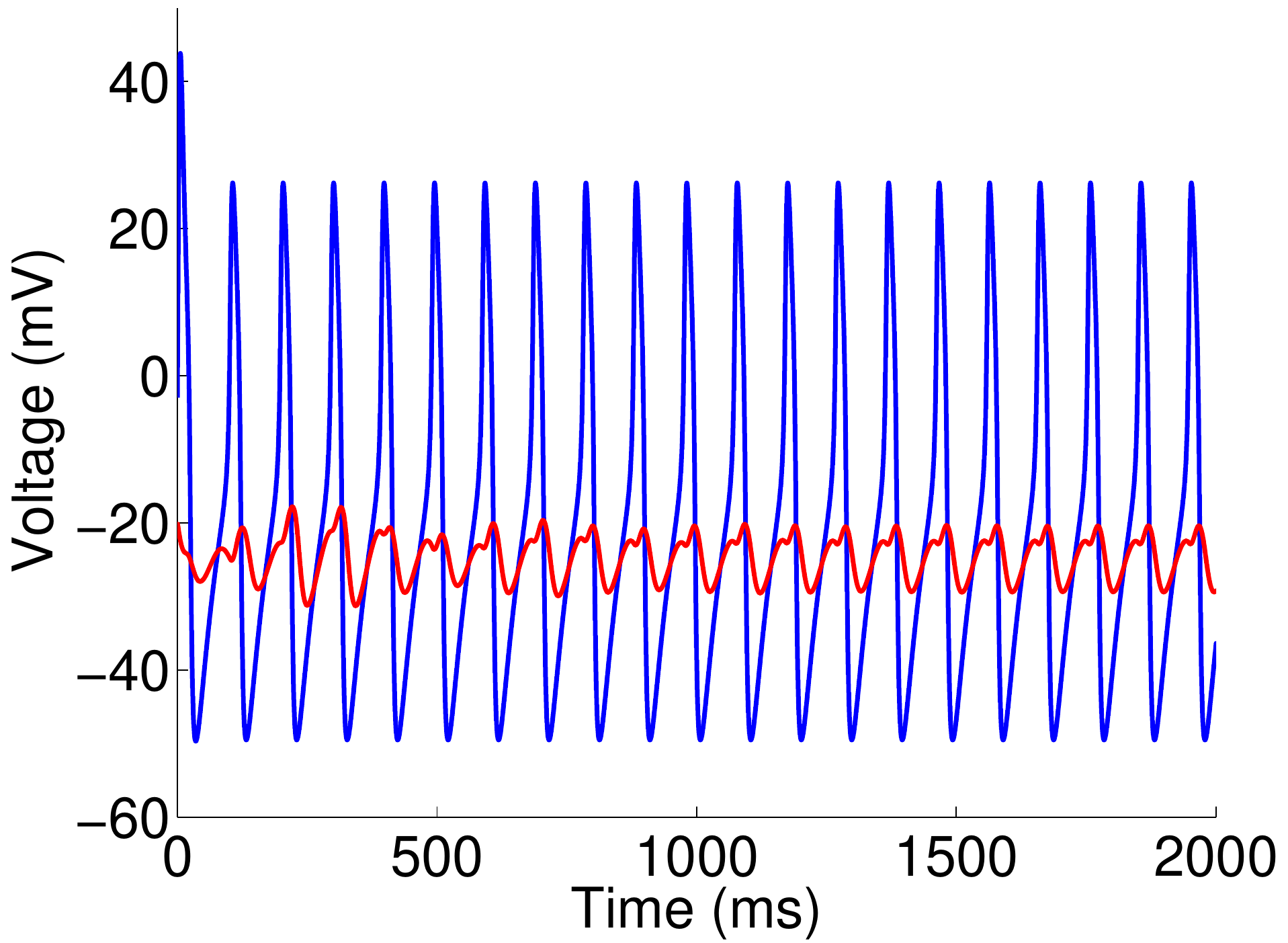} \label{AAS}\\
c) & & & \\
& \includegraphics[scale=.228]{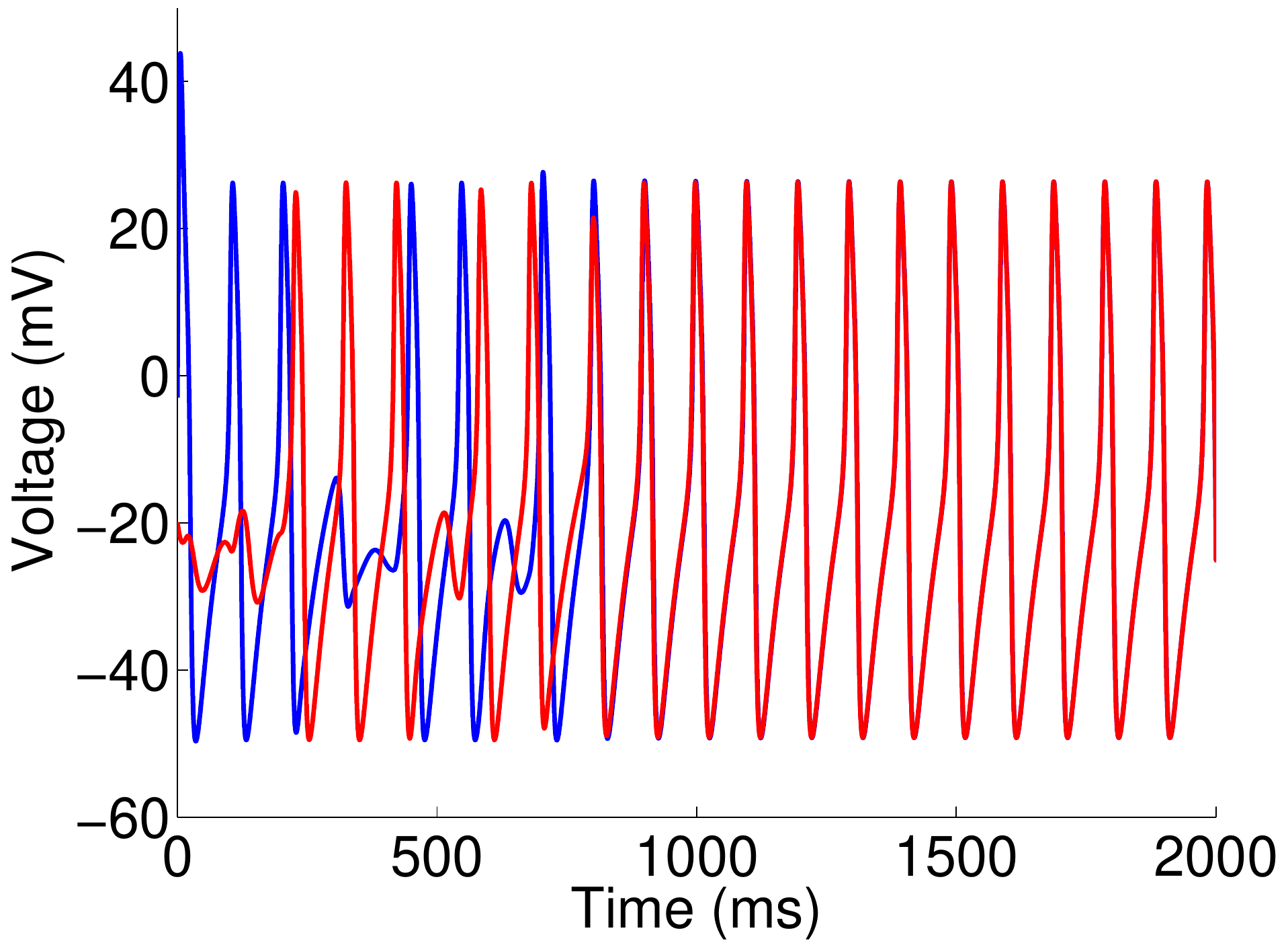}\label{EAS} &&
\end{tabular}
\titlecaption[Stable dynamics of coupled deterministic ML model]{Voltage of neuron 1 is shown in blue and of neuron 2 is shown in red.  In (a), SS dynamics develop after an initial transient.  Parameters are same as in Appendix except that $\Iapp=95.5$ and $\gsyn=0.15$. Anti-phasic AAS arises in (b), such that neuron 1 is experiencing SAO while neuron 2 is experiencing LAO. Parameters are same as in Appendix except that $\Iapp=97.5$ and $\gsyn=0.15$ and the initial conditions $\{v_{1},v_{2},w_{1},w_{2},s_{1},s_{2}\}=$\{-3,-20,0,0.17,0,0\}.  EAS emerges in (c) where voltages of neuron 1 and neuron 2 fire with equal amplitude.  Parameters and initial conditions are same as for b) except $\gsyn=0.2$. Here the voltages gradually become entrained in-phase.}
\label{S}
\end{center}
\end{figure}

When two ML neurons are coupled, there are in addition two different types of periodic behavior. The system is in an asymmetric amplitude state (AAS) when the two voltages are oscillating with different amplitudes. In particular, one neuron experiences large amplitude oscillation (LAO) while the other experiences small amplitude oscillation (SAO). Alternatively, an equal amplitude state (EAS) occurs when both neurons are oscillating with the same amplitude of voltage. Examples of all these states can be found in Fig.~\ref{S}.

If the dynamics of a given data set can be characterized (e.g. AAS, EAS, or SS), then a corresponding parameter regime can be embedded within the MCMC proposal distribution and implemented by simple rejection sampling. In this study, we refine the parameter regime purely based on whether the system is in an oscillatory state, regardless of whether it is AAS or EAS.  Bifurcation analyses were performed computationally using the AUTO library via the software XPPAUT~\citep{xppautbook}.  Evaluations of the model were performed by XPPAUT as needed directly from within MATLAB via the functionality provided by the {\sffamily{xppauttools}} package of the Snifflib JAVA library~\url{sourceforge.net/projects/snifflib}.  All other aspects of analysis were implemented in MATLAB.

\section{MCMC Estimation Method}

\subsection{Metropolis-Hastings Sampling}

MCMC is an iterative random walk method that avoids the need for a closed form of the posterior distribution function.  Rather than a direct method which would require evaluation of a high dimensional integral of the posterior over the parameter space, MCMC behaves much differently, instead returning a sample of parameters from the posterior.  Since these parameter samples are drawn from the desired posterior, consistent estimates of moments such as means and variances of the parameters may be easily calculated.  The Metropolis-Hastings sampler was used which \emph{in lieu} of an exact posterior substitutes a proposal distribution.    We next discuss the proposal, likelihood, and prior distributions supplied code to the MCMC routine. 

The first of these is the proposal, which in Metropolis-Hastings sampling supplies new candidate parameter sets based on the current set in the MCMC chain.  Many of the ML parameters are constrained to be non-negative and in particular we assume this for \gsyn and \Iapp.  However,  we accomplish unconstrained parameter estimation over this support by way of log transformation.  To start, the program is given a set of initial guesses for the log transformed parameters \Iapp and \gsyn. This set of parameters is labeled as ${\theta}_0$. A new candidate parameter set is drawn from a Gaussian distribution with mean centered on these log transformed parameters, ${\theta}_*$.  The standard deviation is a fixed and predetermined value typically referred to as a ``mixing value.''  Then, if the quantity resulting from multiplying the model likelihood by the model prior for the inverse transformed ${\theta}_*$ is greater than that obtained for ${\theta}_0$, then the candidate parameter set is ``accepted'' as the new ${\theta}_0$. Otherwise, the new set ${\theta}_*$, is ``rejected'' and ${\theta}_0$ remains the initial guess. A new parameter set candidate is drawn at random from the same distribution, continuing this cycle until a ${\theta}_*$ is ``accepted.'' The method is continued until after a burn-in period and the chain is determined to have settled into an steady state equilibrium.

\subsection{Construction of Conditioned Likelihood}

Now we discuss implementation of the likelihood.  Analytic and/or efficient forms for the general voltage likelihood for coupled ML are intractable and unavailable. Instead, we construct conditional likelihoods based on feature statistics~\citep{fraser:local1964,fraser2004ancillaries,ghosh:ancillary}.  The first feature statistic we consider is the cumulative power.  

\subsubsection{Cumulative Power}
For a twice differentiable signal $m(t)$, cumulative power $P(t)$ is defined as
\begin{align*}
	P(t) = \int_{0}^{t} \left(m^{\prime\prime}(s)\right)^{2} ds
\end{align*} 
for real finite $t$.  $P(t)$ is a useful feature statistic for a variety of models. For example, square-integrable functions, including bump functions, maxima/minima curves, and saturation curves all have $P(t)\in{\Oh{1}}$.  In contrast, any finite sum of sines and cosines have $P(t)\in{\Ot}$~\citep{quinn}. By way of Fourier representation, this $\Ot$ behavior characterizes a broad array of continuous periodic functions (see proof in Appendix).  Surprisingly, it has even been shown that $P(t)\in{\Ot}$ in a class of linear stochastic dynamical systems lacking differentiability at countably many times~\citep{Bates2012}.  In all such cases, having $P(t)\in{\Ot}$ imparts the intuitive notion that the dynamical system  exhibits a consistent (stationary) duty cycle and accumulates power (on average) at a constant rate.

\begin{figure}[H]
	\centering
	\includegraphics[scale=0.6]{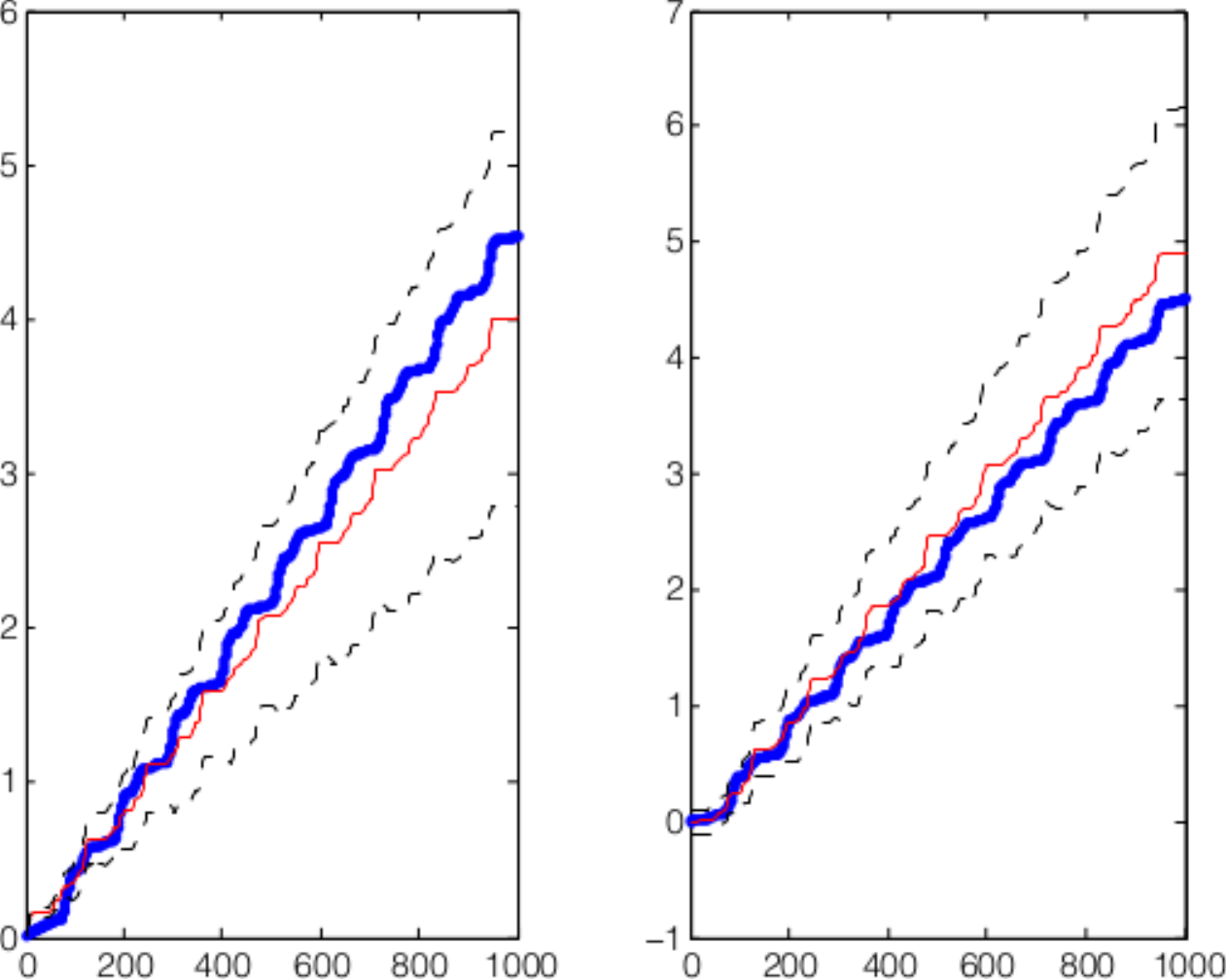}
	\titlecaption[Cumulative power of voltage]{The blue line is the estimated cumulative power of the observed voltage ($P_{\text{data}}(t)$). The red line is the cumulative power of the voltage predicted from the ML model ($P_{\text{model}}(t)$) based on candidate parameter estimates sampled from the proposal distribution. Dashed lines give the $\pm{1}$ standard deviation wedge interval. See Eqn.\ref{LIKELIHOOD} for details.}
	\label{PtofML}
\end{figure}

For periodic oscillating functions such as the voltage in the ML model, the graph of the cumulative power is similar to the red or blue line in Fig.~\ref{PtofML}. 
In Fig.~\ref{PtofML}, the blue line is the estimated cumulative power of the data set $P_{\text{data}}(t)$. The red line is the cumulative power of the model $P_{\text{model}}(t)$ based on parameter estimates with the dashed lines giving the $\pm{1}$ standard deviation interval.(see Eqn.\ref{LIKELIHOOD} for details).  Intuitively, matching  $P_{\text{model}}(t)$ to $P_{\text{data}}(t)$ imparts qualitative matching in terms of amplitude and frequency and may reduce bias, especially in those parameters of the model which determine these features. \Iapp and \gsyn strongly affect both of these features. 

\begin{figure}
\begin{center}
\scalebox{0.4}{\includegraphics{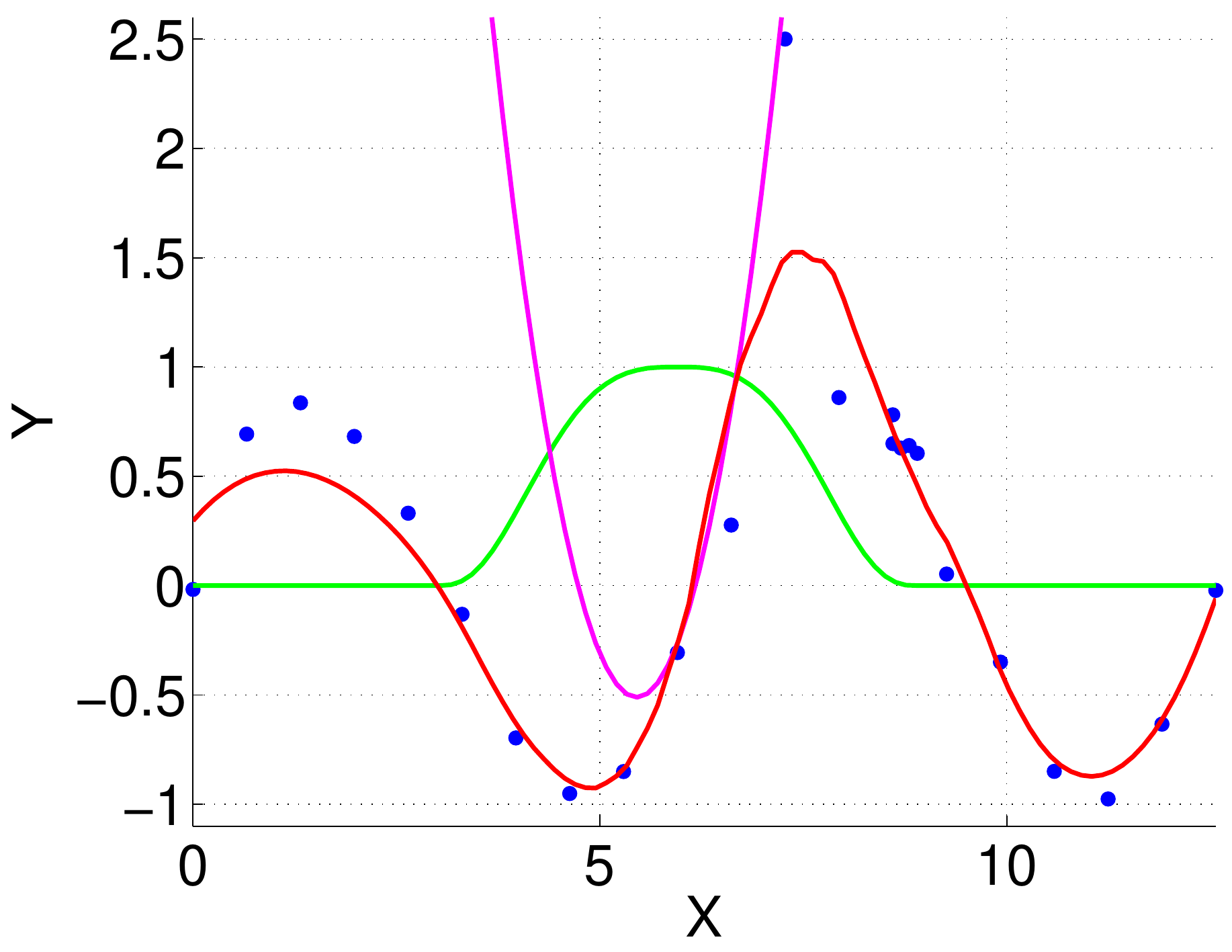}}

\begin{pspicture}(0,0)
	\rput(-4.2,4.7){\begin{minipage}{11ex}\textcolor{black}{Smooth at boundaries}\end{minipage}}
\end{pspicture}

\begin{pspicture}(0,0)
	\rput(3,5.7){\textcolor{black}{Robust to outliers}}
\end{pspicture}

\begin{pspicture}(0,0)
	\rput(4,4.7){\begin{minipage}{10ex}\textcolor{black}{Adapts to sampling density}\end{minipage}}
\end{pspicture}
\vspace{-10ex}
\titlecaption[Example LWPR on noisy unevenly spaced data]{Noisy data (blue) $y_{i} = sin(x_{i}) + \varepsilon_{i}$ where $\varepsilon_{i}\sim{Gaussian(0,0.1)}$.  Locally quadratic regression using a tricube weighting kernel and a nearest-neighbor bandwidth of 53\% (red) goes smoothly through the data coordinates.  The values of the tricube weighting kernel (green) give most weight to the points near the point to be predicted (weights for x=5.95 being shown).  The locally quadratic regression (magenta) may easily be differentiated to generate robust 1st and 2nd derivative estimates at that point.}
\label{LWPR}
\end{center}
\end{figure}

The curves were determined by locally weighted polynomial regression (LWPR) which is a time-domain  smoothing method competitive with frequency domain methods such as Butterworth filtering~\cite{watkin} and wavelets~\citep{fan:local} and are especially convenient when the prediction of time derivatives is desired.  In standard least squares, the over-determined system
\begin{align*}
	Y_{n\times{1}}=X_{n\times{p}}b_{p\times{1}}
\end{align*}	
is taken to have the solution
\begin{align*}	\Hat{b}_{p\times{1}}=\left(X^{\top}_{n\times{p}}X_{n\times{p}}\right)^{-1}X^{\top}_{n\times{p}}Y_{n\times{1}}.
\end{align*}
In contrast, in weighted least squares, each equation is weighted by a diagonal matrix $W_{n\times{n}}$
\begin{align*}
	W_{n\times{n}}Y_{n\times{1}}=W_{n\times{n}}X_{n\times{p}}b_{p\times{1}},
\end{align*}	
which is taken to have the solution
\begin{align*}	\Hat{b}_{p\times{1}}=\left(X^{\top}_{n\times{p}}W_{n\times{n}}X_{n\times{p}}\right)^{-1}X^{\top}_{n\times{p}}W_{n\times{n}}Y_{n\times{1}}.
\end{align*}
The $j$th diagonal element of $W_{n\times{n}}$ is the evaluation of a weighting kernel having a compact support taken at the distance (in time) between the point being predicted and the $j\text{th}$ datum.  Weighting kernels are parameterized by a bandwidth or tuning parameter $h$ which determines the distance beyond which the weights go to zero.  ~\citet{fan:local} gives a comprehensive review of typical weighting kernels used in LWPR. In our application, we choose a tricube weighting kernel and a nearest-neighbor bandwidth which is illustrated in Fig.~\ref{LWPR}.

To determine $P(t)$ from LWPR for each datum  $j\in\left\{1,\ldots,n\right\}$ and for each voltage  $k\in\left\{1,2\right\}$ we calculate,\\
\begin{align}
	\Hat{P}_{\text{data},k}(t_{j})& = \sum_{i=0}^{j}\left[\Delta{t}_{k,i}\left(\Hat{Y}_{\text{data},k}^{\prime\prime}(\given{t_{k,i}}{h_{\text{data}}})\right)^{2}\right]\\
	\Hat{P}_{\text{model},k}(t_{j})& = \Hat{\text{pleak}}_{k}\sum_{i=0}^{j}\left[\Delta{t}_{k,i}\left(\Hat{Y}_{\text{model},k}^{\prime\prime}(\given{t_{k,i}}{h_{\text{model}}})\right)^{2}\right]\\
	\mathcal{L}^{\text{pow}}_{k,j}& = \text{Gaussian}\left(\given{\Hat{P}_{\text{data}}(t_{j})-\Hat{P}_{\text{model}}(t_{j})}{0, \Hat{\text{p0}}_{k} + \Hat{\text{pscale}}_{k}\Hat{P}_{\text{data}}(t_{j})}\right)\label{LIKELIHOOD}
\end{align}
where $\text{Gaussian}(\given{x}{\mu,\sigma})$ is a Gaussian PDF with mean $\mu$ and standard deviation $\sigma$ evaluated at $x$. The second derivative estimates are easily obtained  via the LWPR (see~\citet{fan:local} for details) on each of the voltages. LWPR bandwidths $h_{\text{data}}$ and $h_{\text{model}}$ were determined by generalized cross-validation and fixed during MCMC.  In the current study, we  did not estimate but rather fixed $\Hat{\text{p0}}_{k}$ at the root mean square error from a linear regression of $\Hat{P}_{\text{data}}(t_{j})$ on $t$ which gave satisfactory results.  By similar reasoning, $\Hat{\text{pscale}}_{k}=1$ is used as a reasonable initial guess but was estimated during MCMC.  Note that in this formulation, there is no requirement that the voltages be observed at the same time nor that the spacings between times are equal.

\subsubsection{Mean of the Voltages}

In addition to conditioning on cumulative power, we conditioned the likelihood on the mean voltages.  Specifically, we calculate,\\

\begin{align}
	\mathcal{L}^{\text{mean}}_{k}& = \text{Gaussian}\left(\given{\Bar{V}_{\text{data},k}-\Bar{V}_{\text{model},k}}{\Hat{\text{mstd}}_{k}}\right)
\end{align}
where $\Bar{V}_{\text{data}}$ is the mean of the observed voltage data for the $k$th neuron,  $\Bar{V}_{\text{model}}$  is the mean of the predicted voltage data for the $k$th neuron, and $\Hat{\text{mstd}}_{k}$ is the estimated standard deviation of the mean voltage residual for the $k$th neuron.  A reasonable initial guess for the $\Hat{\text{mstd}}_{k}$ is 10mV.  The conditioned model likelihood is thus obtained by multiplication of the $\mathcal{L}^{\text{pow}}_{k}$ and  $\mathcal{L}^{\text{mean}}_{k}$ likelihoods,
\begin{align}
	\mathcal{L}(\theta_{*})& = \Pi_{k=1}^{2}\left(\Pi_{j=1}^{n}\mathcal{L}^{\text{pow}}_{k,j}\right)\cdot{\mathcal{L}^{\text{mean}}_{k}}.
\label{modellike}
\end{align}
where $\theta_{*}=\left\{\Hat{\Iapp},\Hat{\gsyn},\Hat{\text{pscale}}_{1},\Hat{\text{pscale}}_{2},\Hat{\text{pleak}}_{1},\Hat{\text{pleak}}_{2}, \Hat{\text{mstd}}_{2}, \Hat{\text{mstd}}_{2}\right\}$.

\subsection{Construction of Parameter Rejection Region}
Bifurcation analysis can be used to learn about the behavior of the system, specifically where the system is EAS/AAS or {SS}. Assuming we know the state of the data set, we can use the bifurcation analysis to limit the parameter combinations over which MCMC tests. This should lead to a more accurate approximation of \Iapp and \gsyn.

 One reason we choose to construct parameter boundaries based on the state of the system is because of the vast difference in cumulative power for EAS/AAS versus {SS} solutions.  Since MCMC calculates probability based on cumulative power, proposed points in a different state than the real data set will always lead to a rejection. In addition, it is easy to observe whether a data set is oscillating without knowing the parameter values.

\subsubsection{Limit on \gsyn}
\begin{figure}[H]
\begin{center}
\begin{tabular}{rcrc}
a) & \ensuremath{\gsyn=4.0} &b) &\ensuremath{\gsyn=7.0}\\
& \includegraphics[scale=.24]{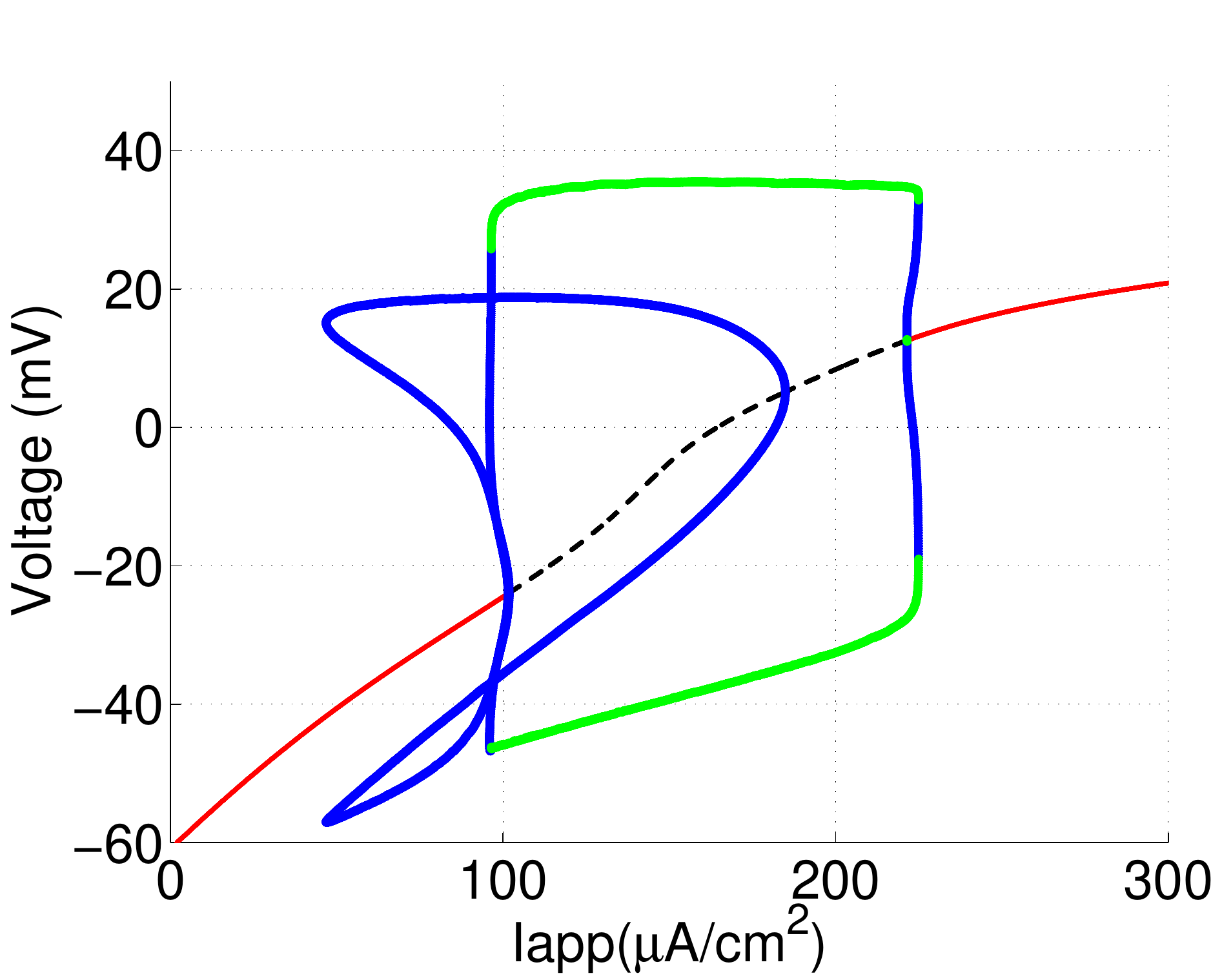} & &
\includegraphics[scale=.24]{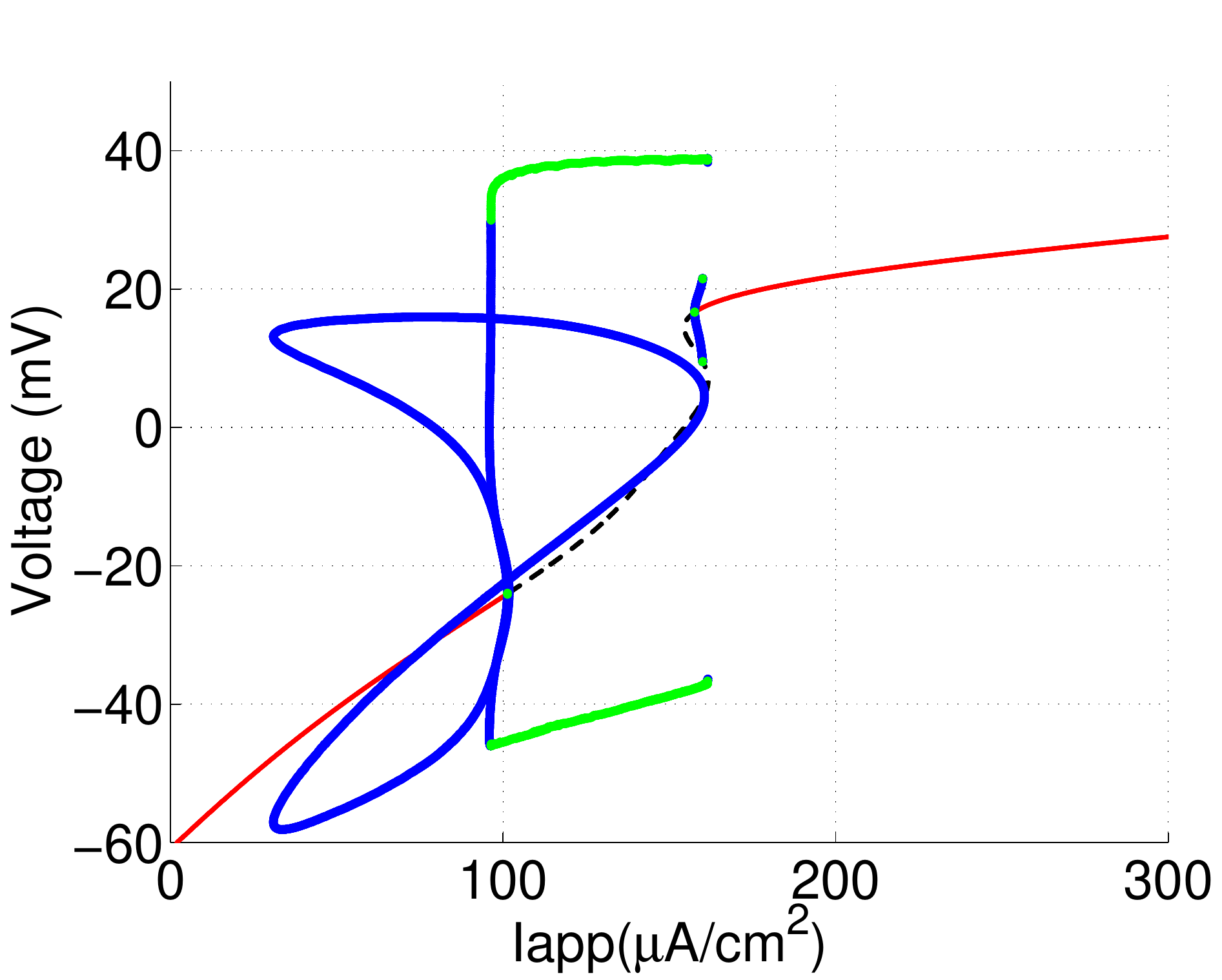}\\
c) &\ensuremath{\gsyn=9.08} & d) &\ensuremath{\gsyn=10.0}\\
& \includegraphics[scale=.24]{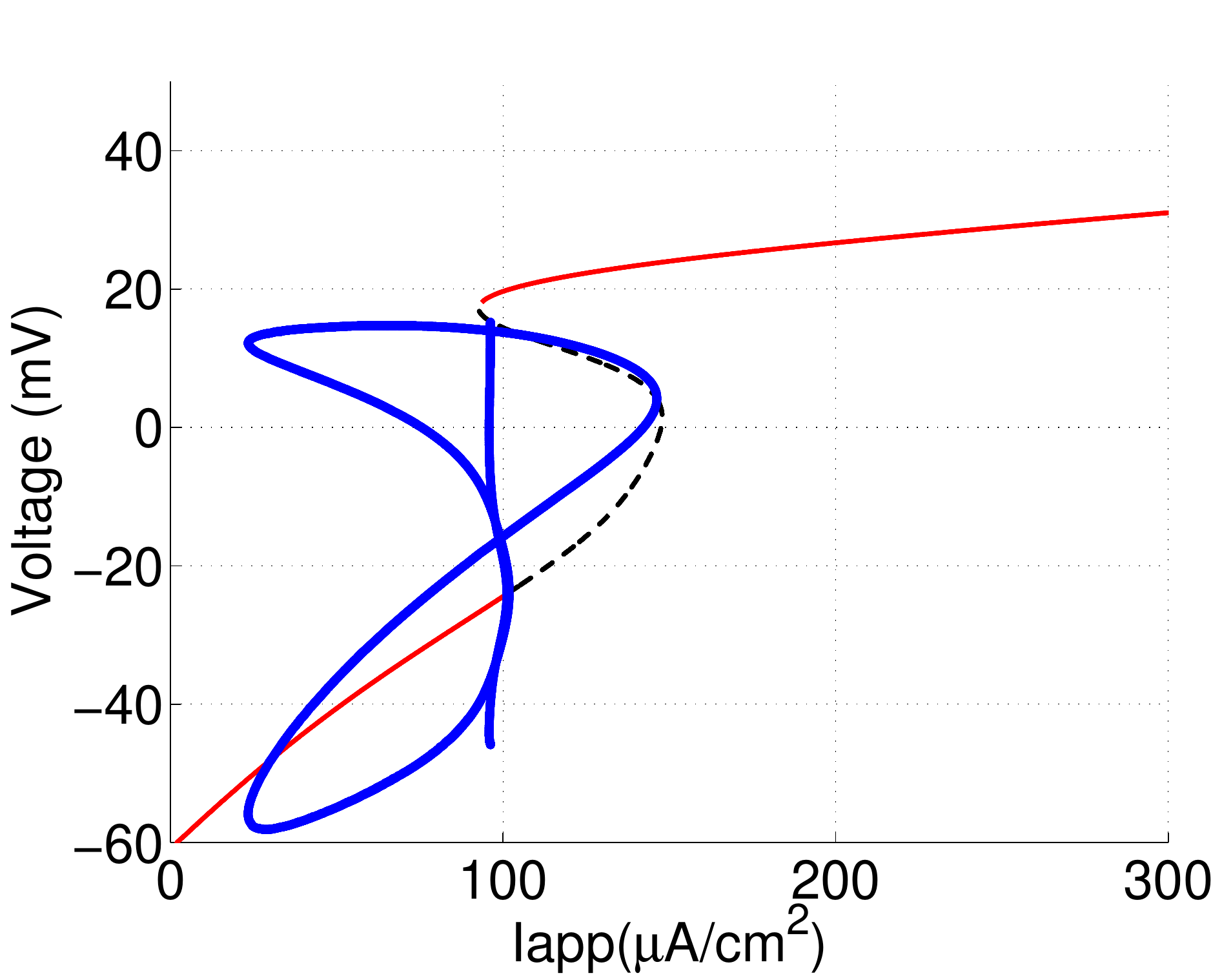} & & \includegraphics[scale=.24]{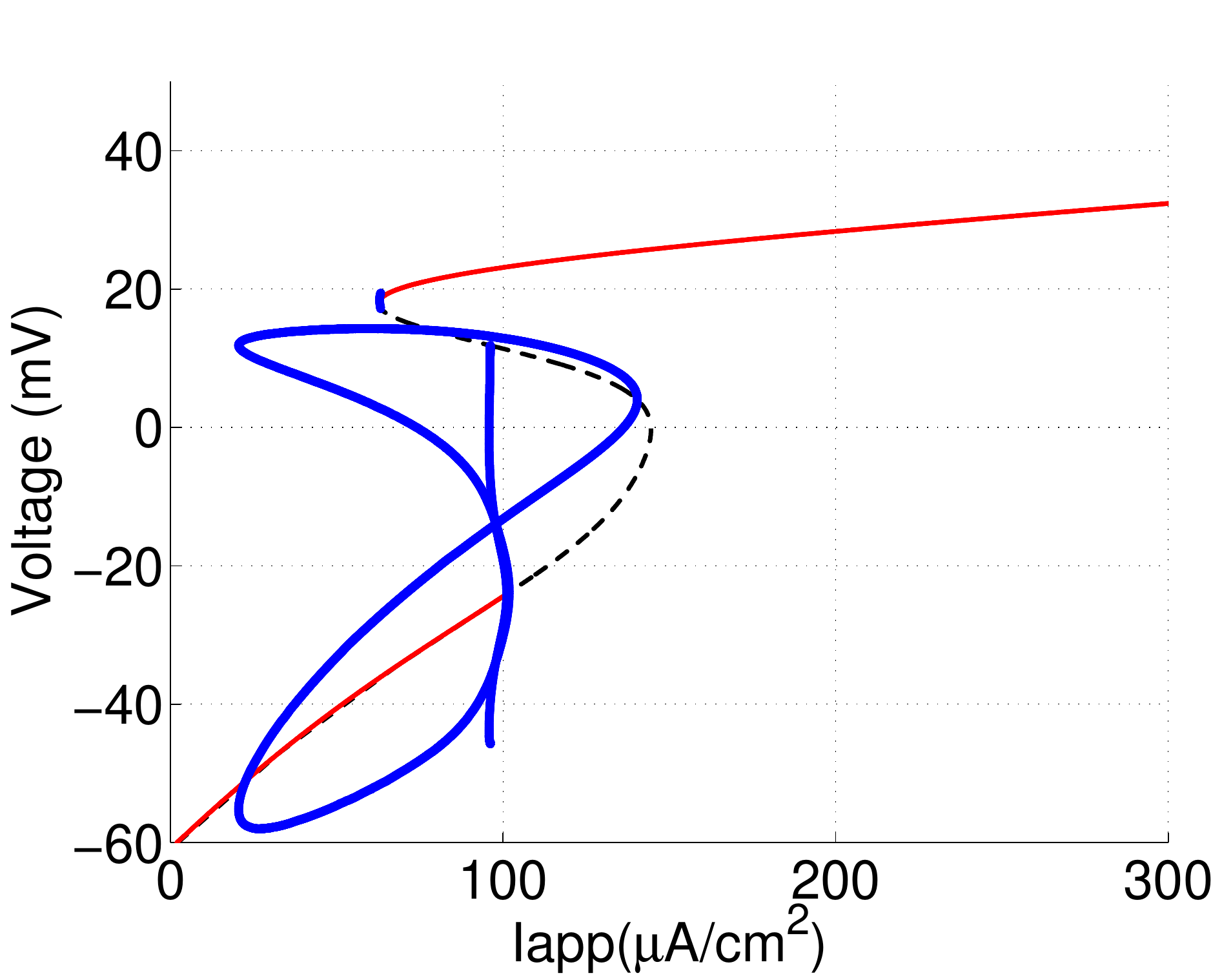}
\end{tabular}
\titlecaption[Loss of stable periodic orbits]{In a) a there is a wide range of \Iapp values for which stable periodic orbits are possible. In b) the range of \Iapp values leading to stable periodic orbits is considerably smaller. In c) stable periodic orbits are not possible and (as seen in d) do not return for larger values of \Iapp.}
\label{B}
\end{center}
\end{figure}

A system will always exhibit the behavior of the stable solution for a set of parameter values. For instance, in Fig.~\ref{B}(a) for an \Iapp value of 150, there is a stable periodic solution, an unstable periodic solution, and an unstable steady state solution present. Depending on the initial conditions used, the system may move briefly toward the unstable solutions, but eventually it will settle to exhibit behavior consistent with the stable periodic solution. Therefore, long-term oscillating solutions will only be present for parameter values corresponding with a stable periodic branch in the bifurcation diagram. Starting at a \gsyn of 0, it can be observed that the stable periodic branch of the system becomes shorter as \gsyn increases. This can be seen in Fig.~\ref{B}. The length of the branch decreases between \gsyn values of 4 and 7 before disappearing altogether at a \gsyn value of 10. Narrowing this further, we find that the value of \gsyn for which the stable periodic branch ceases to exist is about 9.08. At this point, the system is always {SS}.  Therefore, based on our observed behavior, oscillatory solutions only occur for \gsyn values less than 9.08. To maintain biological realism, we also impose $\gsyn\geq{0}$ which provides a lower bound.   

\subsubsection{Limit on \Iapp}
With an established range for \gsyn, we wish to constrain the corresponding \Iapp values for EAS/AAS oscillations. To do this, we carried out two-parameter bifurcation continuation of the limit point (LP) bifurcations.  This revealed an envelope enclosing a similar two-parameter  Hopf continuation. Hopf bifurcations can be of either the  sub-critical  or super-critical type with stable periodic orbits possibly emerging below or above the Hopf bifurcation.  We observe the LP continuation envelope extending below the Hopf envelope for smaller values of \Iapp and extends above the Hopf envelope for larger values of \Iapp.   Outside of the LP region, SS is the only stable behavior. Table~\ref{POLYREGION} shows the final oscillation parameter ranges found. It can be observed that as \gsyn increases, the LP envelopes draw closer until eventually intersecting around \gsyn=9.08 and \Iapp=95.840. This is consistent with our \gsyn boundary. 

 \begin{table}
\centering
\begin{tabular}{|c|c||c|c|}
\hline 
\textbf{$\mathbf{\gsyn}$} & \textbf{$\mathbf{\Iapp}$} & \textbf{$\mathbf{\gsyn}$} & \textbf{$\mathbf{\Iapp}$}\\ \hline \hline
 0.007 & 238.382 &  9.019 & 95.842\\
 0.832 & 238.097 &  8.199 & 95.867\\
 1.657 & 237.885 &  7.380 & 95.889\\
 2.482 & 236.004 &  6.561 & 95.909\\
 3.306 & 231.454 &  5.742 & 95.925\\
 4.131 & 223.402 &  4.923 & 95.938\\
 4.956 & 211.055 &  4.103 & 95.945\\
 5.781 & 193.919 &  3.284 & 95.946\\
 6.606 & 172.617 &  2.465 & 95.939\\
 7.430 & 148.735 &  1.646 & 95.918\\
 8.255 & 123.144 &  0.826 & 95.871\\
 9.080 & 95.840 &   0.000 & 95.724\\\hline
\end{tabular}
\titlecaption[Proposal rejection look-up table]{These 24 values give good approximation to the region where EAS/AAS oscillations are possible (see Fig.~\ref{ROLL}). Defining a closed region, these few values allow computationally efficient rejection sampling within the proposal distribution function.}
\label{POLYREGION} 
\end{table}

\begin{center}
\begin{figure}[H]
\centering
\includegraphics[scale=0.7,viewport=1.0in 2.5in 7.5in 8.25in]{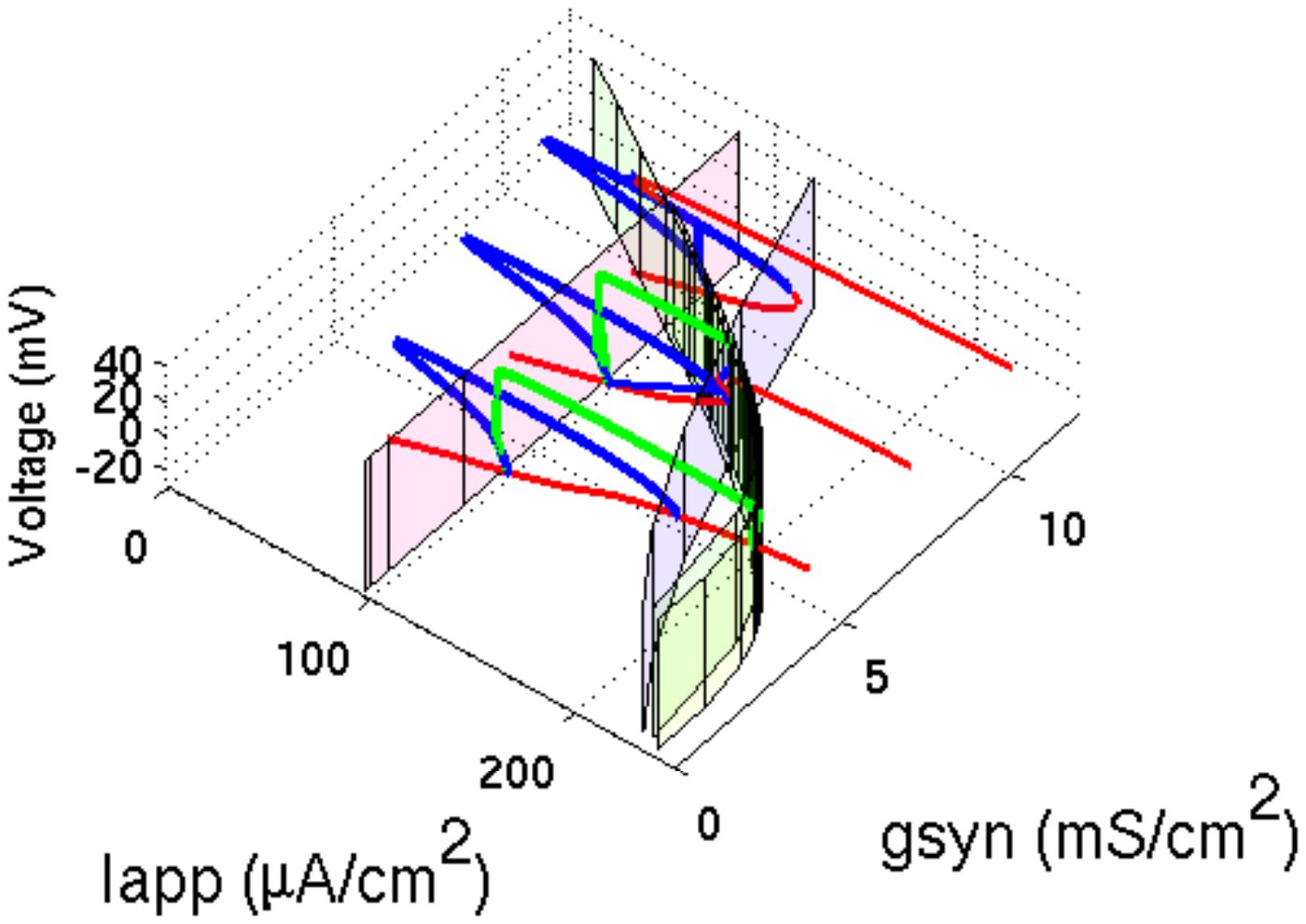}\\
\titlecaption[3D bifurcation diagram]{Hopf and LP enclosures (see 2D projection in Fig.~\ref{ROLL}) are shown intersected by planar bifurcation diagrams taken at cross-sections $g_{syn}=$ 4, 7, and 10.}
\label{param2}
\end{figure}
\end{center}

\subsubsection{Rejection Sampling}

\begin{figure}[H]
\begin{center}
\includegraphics[scale=0.6]{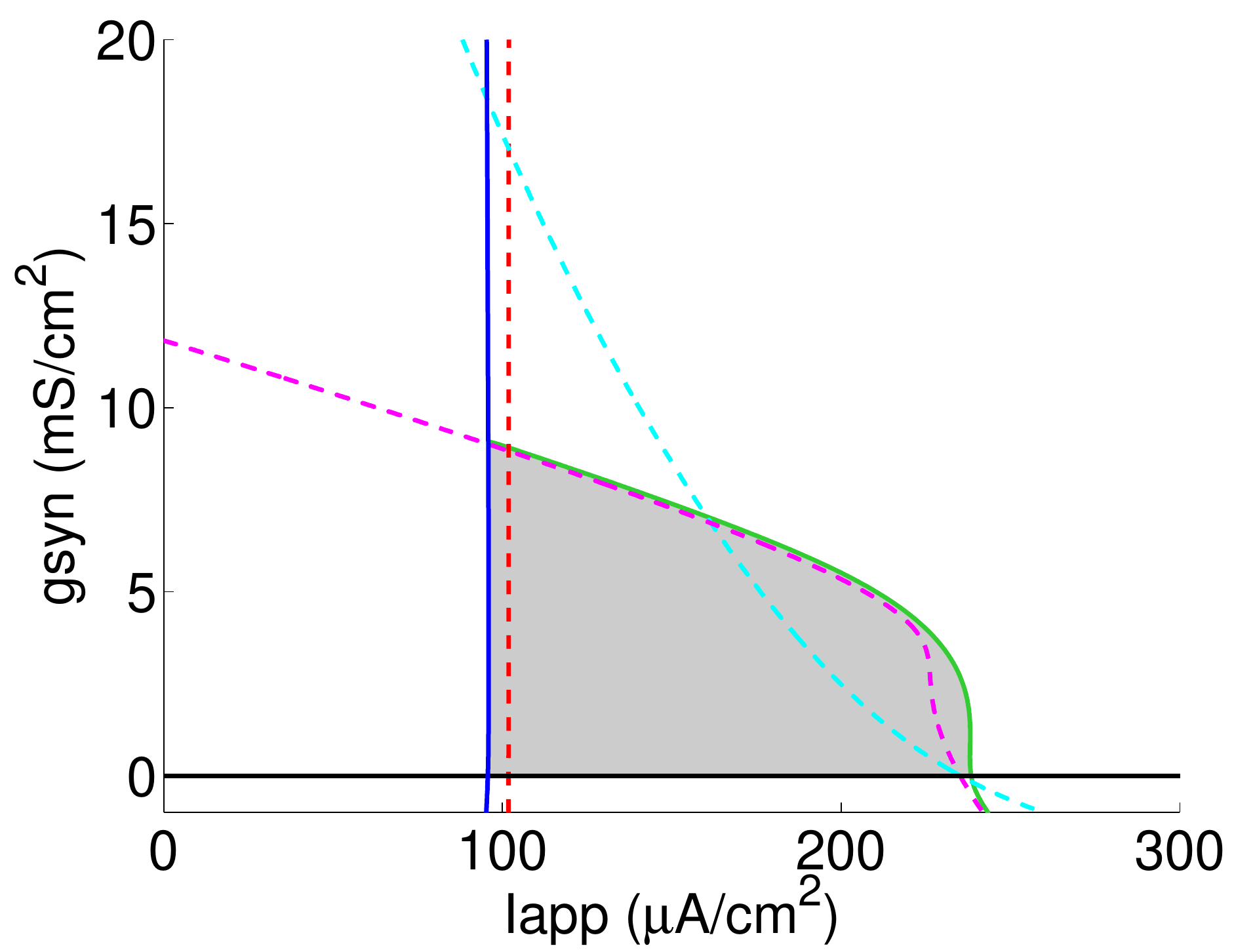}%
\begin{pspicture}(0,0)%
\rput(-5.0,2.5){\begin{minipage}{18ex}\centering NO ROLLBACK\\EAS/AAS\\$P(t)\in\Ot$\end{minipage}}
\rput(-7.5,6){\rotatebox{45}{\begin{minipage}{20ex}\centering ROLLBACK\\\Large{SS $\quad P(t)\in\Oh{1}$}\end{minipage}}}
\end{pspicture}%
\titlecaption[Rejection region for oscillating solutions (AAS and EAS)]{As in Fig.~\ref{B}a), the Hopf bifurcation continuations (cyan, red, and magenta dash) exhibit sub-critical (red) or super-critical (magenta) stable periodic orbits.  Accordingly, stable periodic orbits are found for any parameter combinations bounded within a slightly larger region (gray) having an area $\approx{933.75}~\mu{A}\cdot{mS}\cdot{cm^{-4}}$ and coinciding with LP bifurcation continuations (green and blue solid).  Outside of this region, SS is the only stable behavior.} 
\label{ROLL}
\end{center}
\end{figure}

\begin{figure}[H]
\begin{center}
\includegraphics[scale=0.6]{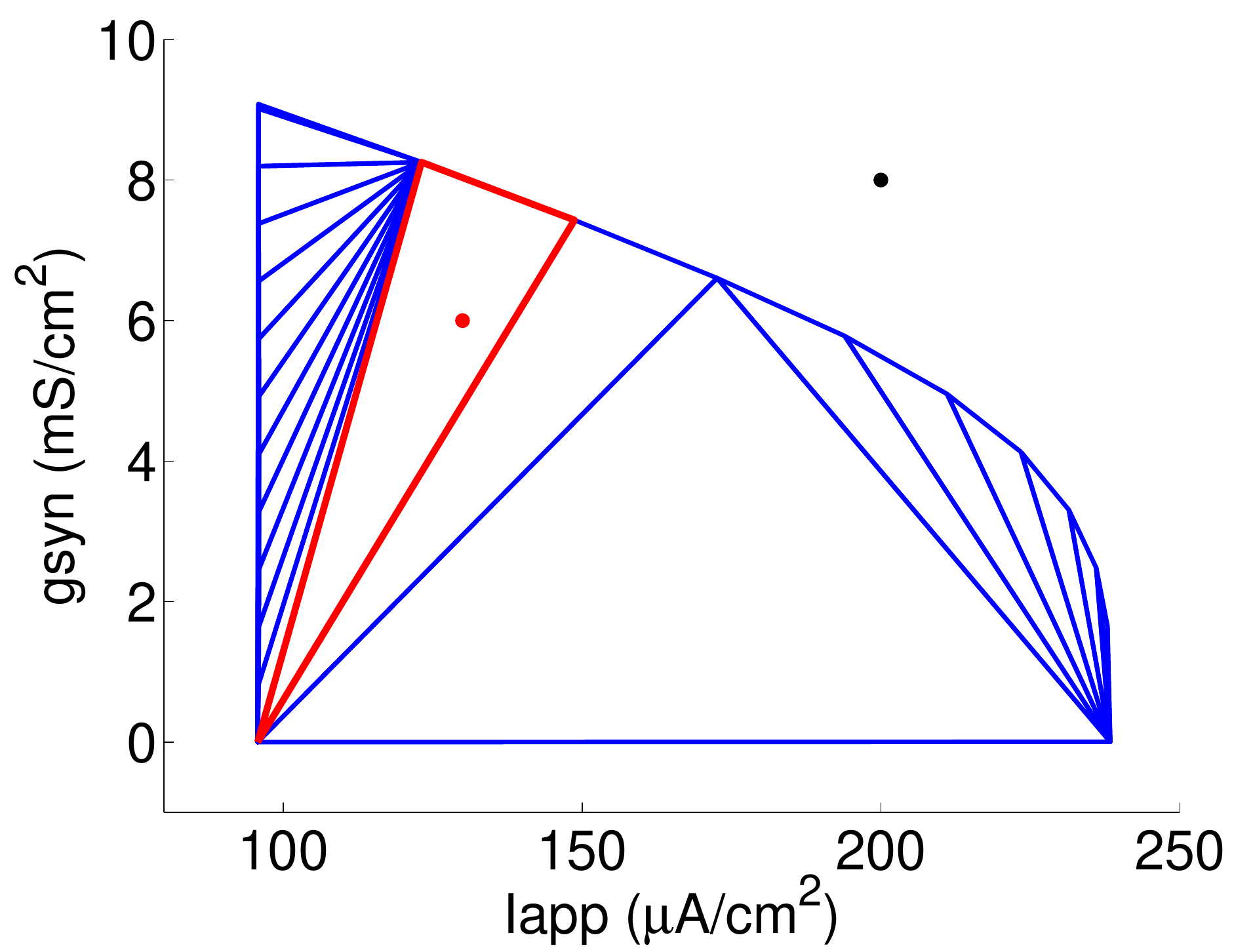}
\titlecaption[Delaunay triangulation of ``NO ROLLBACK'' region]{In an approach which extends naturally to higher parameter dimensions, the LP region in Fig.~\ref{ROLL} is decomposed into triangular regions (blue). An acceptable candidate parameter set (red) is efficiently identifiable by being interior to one of the triangles. In contrast, an unacceptable parameter set (black) is not interior to any of the triangles and would result in a rollback.}
\label{DELAUNAYTRI}
\end{center}
\end{figure}

Determining if a parameter candidate is within the 2D region specified in Table~\ref{POLYREGION} can be easily accomplished by a variety of methods.  However, we have interest to extending our methods to more parameters than just \Iapp and \gsyn.  Therefore Delaunay triangulation was adopted as a general approach to the in/out region testing. In this approach the possibly D-dimensional region is decomposed into simplices (triangles in 2D) having $D+1$ vertices  such that the $D$ dimensional circumspheres about any vertex do not contain any other points~\citep{cignoni_dewall}.  The number of simplices is determined by the number of points defining the region boundaries which is our case result from parameter bifurcation continuation via XPPAUT.  This triangulation process is efficient and may even be performed offline prior to MCMC and thus represents negligible computational expense. As seen in Fig.~\ref{DELAUNAYTRI}, a  parameter candidate generated by a random draw from the proposal distribution is determined to be within this region if it is internal to any one of this simplices.  For a region determined by $m$ continuation points, this is an $\Oh{m}$ operation and quite efficient.   

\subsection{Elicitation of Prior}

A valid posterior probability requires multiplication of the likelihood in Eqn.~\ref{modellike} by the prior probability of the candidate parameters.  In determining an informative prior for \Iapp and \gsyn, the nonlinear nature of the model should not be underestimated.  Typically parameters interact in highly nonlinear ways to establish the behavior of the system.  The different regions of parameter space can, in principle, be assigned distinct prior probabilities.  Another approach would be to define a density that varies smoothly and continuously over the parameter support.  However, in either case the $p$-dimensional improper integral, taken over the support of all parameters, of the joint prior distribution must evaluate to 1.0.  If some regimes, but not others, are to be considered then an p-dimensional density over a typically irregular region must be defined.  There is typically little guidance about the shape of such regions, in particular such regions may not even be  simply connected.  

Our approach regarding elicitation of the prior follows from the typical goals of mathematical modeler who wants to know if their model is capable of predicting the response while capturing certain key \emph{features} of interest.  Consider that a constant voltage set at the average voltage of a neuron firing action potentials may in fact provide a relatively high predictive accuracy since it corresponds to the least squares linear regression.  However the average voltage model, clearly would be unacceptable because it predicts poorly the additional \emph{features} of amplitude and frequency.  We assume researchers will have sufficient expertise to know that their data exhibits action-potentials and feel justified in setting to zero the prior probability of any parameters which do not predict such behavior.  

We find that this ``NO ROLLBACK'' region is a simply connected finite region  (having area of $A\approx{933.75}~\mu{A}\cdot{mS}\cdot{cm^{-4}}$). We assume flat uniform prior over this EAS/AAS region.  This is the appropriate approach when an investigator is interested in inference for a model conditioned on a prior belief that the neurons are indeed firing action potentials, but has no further prior information about values of  \Iapp and \gsyn.  Such flat priors are commonly called uninformative priors and are common choice of prior in Bayesian statistical analysis.  The ``no rollback'' region of Fig.~\ref{ROLL} is assumed to have prior probability equal to $1/A$ while all other parameter values have prior probability zero\footnote{Any positive constant (e.g. 1.0) would serve equally well.}.  Non-uniform priors over this region could be defined which would have the effect to drive parameter estimates towards low/high values of \Iapp and \gsyn. Non-uniform priors are beyond the scope of the present paper but are areas for continuing research by our group.   

At each iteration of MCMC, candidate parameters generated by the proposal distribution were tested to determine if they resided within the ``NO ROLLBACK'' region by using the MATLAB function \texttt{pointLocation} with the look-up values in Table~\ref{POLYREGION}.  This approach was highly efficient and did not contribute significantly to the computational effort.  Candidate parameters falling outside the region could be immediately rejected avoiding any further computational expense on that iteration.  

\subsection{Proposal Distribution}

Metropolis-Hastings sampling in MCMC is distinguished by allowing for the next candidate parameters in the chain to be sampled from a proposal distribution.  This is convenient in applications such as ours where an analytic form for the posterior is either unknown, intractable, or computationally inefficient.  Typically a proposal distribution is selected which i) shares the same support as the parameter(s) of interest and ii) is easy to evaluate.  Typically so-called location-scale distributions are utilized so that the next guess is drawn from a distribution may be \emph{centered} on the most recent element in the chain.  The \emph{scale} of the proposal is often called a \emph{mixing} value as it controls, in part, how vigorously parameter space is traversed.  The \emph{mixing} values may also be adaptively tuned during MCMC.  This so-called adaptive MCMC is especially useful or even required in nonlinear problems where characteristically the sensitivity of the likelihood and/or prior to small changes in parameter values can change drastically as the parameter space is traversed.

Both \Iapp and \gsyn are physiological parameters.  In our application we were only interested in excitatory input and positive synaptic conductances.  To achieve unconstrained sampling of the positive semi-infinite support of these parameters we estimated the log of the parameters' values. The proposal distribution for both log-parameter was assumed to be a Gaussian.  Mixing values (standard deviations) for these proposals was initially assigned to 0.01, but then adaptively varied by the operator down to 0.001 after burn-in was determined to have occurred.

\subsection{MCMC Without Rejection Sampling}

\begin{center}
\begin{figure}[H]
\centering
\includegraphics[scale=0.6]{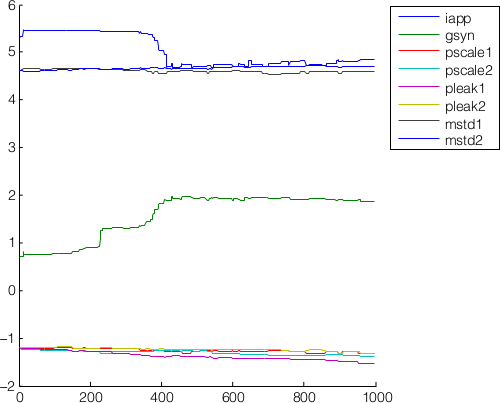}\\
\titlecaption[MCMC trace plot]{The blue line is the natural log of the current estimation of the \Iapp. The green line is the natural log of the current estimation of the \gsyn. The burn-in phase of this particular trace, where the current estimate takes larger steps to find a more likely estimate, ends after 450 iterations. The post-burn-in phase is after this point, where the program focuses in on a small ranges of highly likely estimates.}
\label{mcmctrace}
\end{figure}
\end{center}

For purposes of comparison, we now illustrate the MCMC estimation of \Iapp and \gsyn without rejection sampling by ignoring completely the feasibility region illustrated in Fig.~\ref{ROLL}.  In Fig.~\ref{mcmctrace}, the parameter trace of the MCMC program is shown for one thousand iterations, where the moving blue line is the natural log of the current \Iapp estimate and the moving green line is the natural log of the current \gsyn estimate. The flat sections where the guess is unchanged for a number of time steps represent a string of rejects where the previous guess, ${\theta}_0$, was kept as the most probable estimate. The regions where the trace jumps to a new estimate is where a new set of estimates was accepted and replaced the previous. 

The trace can be split into two sections. First is the burn-in phase, where the current estimate is jumping relatively far from one set of parameter guesses to another as it comes closer to the most likely set of values. The second part is the post-burn-in phase, when the estimate limits itself to a very small range and remains more or less stationary around the mostly likely sets of values. For Fig.~\ref{mcmctrace}, the burn-in phase ends at about 450 time steps. To find the best estimate, we take the average of the post-burn-in phase after this point. We ran the MCMC routine for five trials in the process outlined above. For known parameters and the listed initial guesses, the program finds the data in Table~\ref{paramtable1}.

\begin{table}[H]
\begin{center}
\begin{tabular}{rlrl}
a) &Data with $I_{app}=120$, Initial guess: 220 & b) & Data with $g_{syn}=7.5$, Initial guess: 1\\
& \begin{tabular}{|c | c | c |}
			\hline
			\textbf{Trial} & \textbf{Mean} & \textbf{Std. Dev.} \\ \hline 					\hline
			1 & 238.30 & 1.00 \\ \hline
				2 & 242.53 & 1.00 \\ \hline
				3 & 226.71  & 1.00  \\ \hline
				4 & 240.80 & 1.00\\ \hline
				5 & 238.11 & 1.00  \\ \hline
			\textbf{Overall} & \textbf{237.29} & \textbf{1.00}  \\ \hline
			\multicolumn{2}{|c|}{\textbf{Percent Error}} & \textbf{97.7\%}  \\ 				\hline
\end{tabular} & & %
\begin{tabular}{|c | c | c |}
			\hline
			\textbf{Trial} & \textbf{Mean} & \textbf{Std. Dev.} \\ \hline 					\hline
			1 & 0.00 & 1.83 \\ \hline
				2 & 1.42 & 1.04\\ \hline
				3 & 3.77 & 1.00 \\ \hline
				4 & 1.53 & 1.03 \\ \hline
				5 & 0.00 & 2.57 \\ \hline
			\textbf{Overall} & \textbf{1.35} & \textbf{1.49} \\ \hline
			\multicolumn{2}{|c|}{\textbf{Percent Error}} & \textbf{82.1\%}  \\ 				\hline
\end{tabular}
\end{tabular}
\titlecaption[\Iapp and \gsyn estimation results without rejection sampling]{Means and standard deviations were easily determined as sample statistics calculated over the MCMC trace after burn-in.}
\label{paramtable1}
\end{center}
\end{table}

Obviously, these estimations are less than satisfactory with high percent errors. For both parameters, the program fails to approach the true parameters or deviate from the initial guesses. In addition to parameter sets producing EAS/AAS, elements of the MCMC chain were also occasionally {SS} parameter sets.  This created flat-line cumulative power and mean voltage graphs. While this deviation from the oscillating state of the system was often rejected by MCMC, it typically lead to variance components ($\text{pscale}_{*},\text{mstd}_{*}$) being  over-estimated and if the range searched for more probable values was too small, the chain often become trapped -- predicting only implausible parameter sets.  By letting the program run without informed parameter constraints , MCMC produces unrealistic parameter estimates for different states.  Thus, to find the feasible region where the EAS/AAS solutions exists, we look to bifurcation analysis.

\subsection{MCMC With Rejection Sampling}
Five more trials of this MCMC code were run in the same manner described earlier, but in combination with  rejection sampling. The resulting estimates can be found in Table~\ref{finalresults}. These trials are much more accurate than those without a limited parameter space, as can be seen from the dramatic decrease in percent error from Table~\ref{paramtable1}. Both the estimations of \Iapp and \gsyn moved farther from the initial guess value and gravitated toward the actual value. These results are especially satisfying because the initial guesses were so far from the actual data values.  

\begin{table}[H]
\begin{center}
\begin{tabular}{rlrl}
a) & Data with $I_{app}=120$, Initial guess: 220& b) &Data with $g_{syn}=7.5$, Initial guess: 1\\

&\begin{tabular}{|c | c | c |}
			\hline
			\textbf{Trial} & \textbf{Mean} & \textbf{Std. Dev.} \\ \hline 					\hline
			1 & 102.08 & 1.03 \\ \hline
			2 & 109.32 & 1.10 \\ \hline
			3 & 98.79 & 1.02 \\ \hline
			4 & 110.67 & 1.07\\ \hline
			5 & 100.23 & 1.02  \\ \hline
			\textbf{Overall} & \textbf{104.22} & \textbf{1.05}  \\ \hline
			\multicolumn{2}{|c|}{\textbf{Percent Error}} & \textbf{13.6\%}  \\ 				\hline
\end{tabular}
&&
\begin{tabular}{|c | c | c |}
			\hline
			\textbf{Trial} & \textbf{Mean} & \textbf{Std. Dev.} \\ \hline 					\hline
			1 & 8.17 & 1.03 \\ \hline
			2 & 7.62 & 1.10\\ \hline
			3 & 8.43 & 1.03 \\ \hline
			4 & 7.74 & 1.05 \\ \hline
			5 & 8.34 & 1.02 \\ \hline
			\textbf{Overall} & \textbf{8.06} & \textbf{1.05} \\ \hline
			\multicolumn{2}{|c|}{\textbf{Percent Error}} & \textbf{7.46\%}  \\ 				\hline
\end{tabular}
\end{tabular}
\titlecaption[$\Hat{\Iapp}$ and $\Hat{\gsyn}$ using proposal with rejection sampling]{Means and standard deviations were easily determined as sample statistics calculated over the MCMC trace after burn-in. Compared to results without rejection sampling (see Table~\ref{paramtable1}), these estimates are less biased.}
\label{finalresults}
\end{center}
\end{table}

\begin{figure}[H]
	\begin{center}
	\includegraphics[scale=0.6]{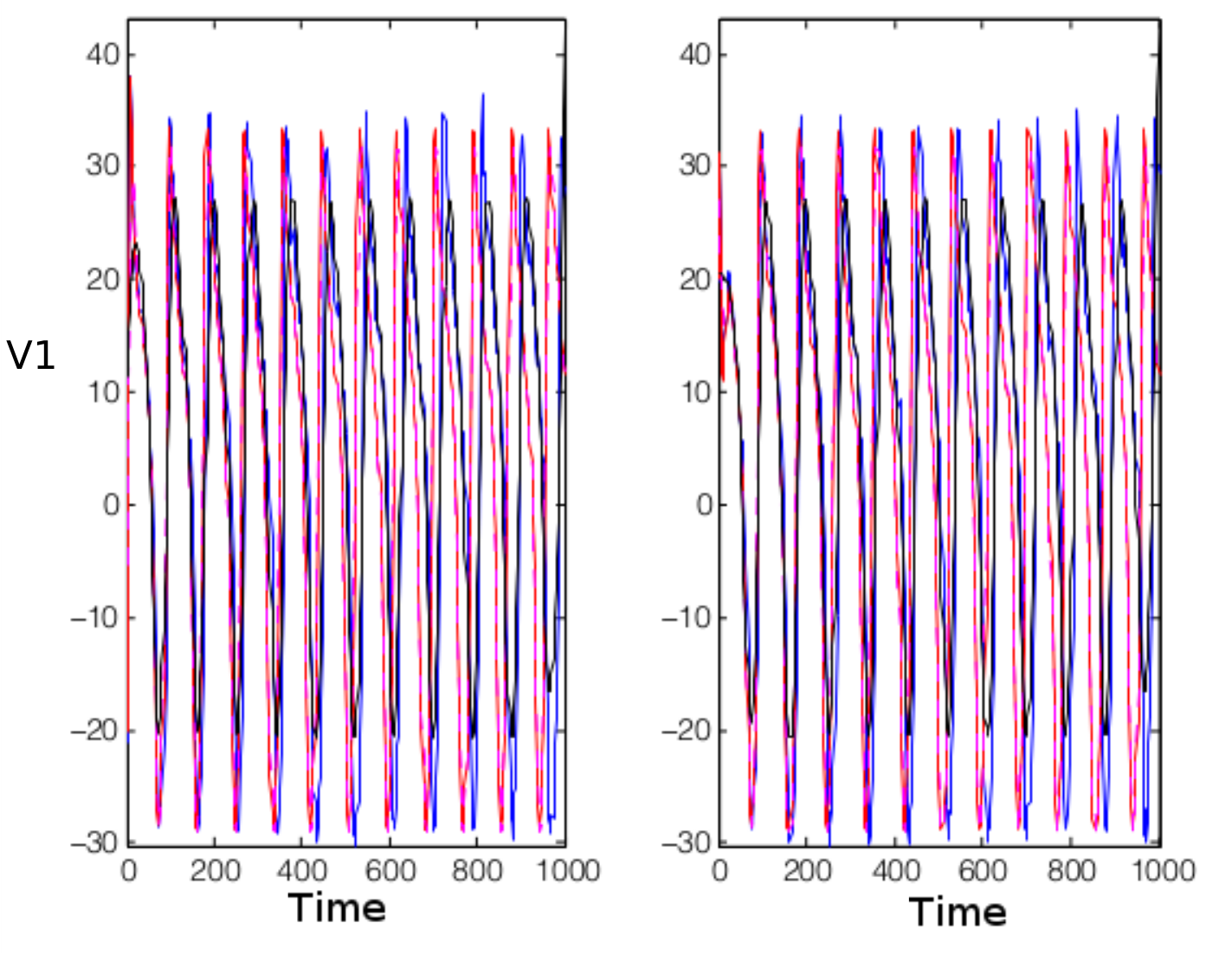}
	\titlecaption[Predicted voltage vs.~time]{ The blue line is the true voltage for the data given, and the red line is the voltage for the current parameter estimation.}
	\label{TrialOne3}
	\end{center}
\end{figure}

\section{Conclusion}

Estimating parameters with MCMC alone is not enough to procure good results for a dynamical system with different behavioral states. One method we find to be effective for the coupled ML model is limiting the parameter region being considered in MCMC using rejection sampling. This can be accomplished through first knowing the state of a data set, then determining parameter boundaries of that state using bifurcation continuation analysis. This removes proposed parameter sets of MCMC that are behaviorally far from the data. In unrestrained MCMC, these sets make it difficult for estimates to ever be ``accepted,'' because the proposed conditioning statistics are far from that of the data. By applying this method to the estimation problem, there is a large decrease in percent error of the estimation parameters, as can be seen in Tables~\ref{paramtable1} and~\ref{finalresults}.  With the primary goal of estimating \Iapp and \gsyn, we introduced the concomittant parameters $\Hat{\text{pscale}}_{1}$, $\Hat{\text{pscale}}_{2}$, $\Hat{\text{pleak}}_{1}$, $\Hat{\text{pleak}}_{2}$, $\Hat{\text{mstd}}_{2}$, $\Hat{\text{mstd}}_{2}$ which linked to the parmaters of interest through specific choices modeling the lack of fit of the predicted voltage and its cumulative power.  As is standard practice in Bayesian latent or strucutral equation modeling, these extra parameters effectively represent an hypothesized latent covariance relation among the parameters of interest and are considered just as integral to model design as the choice of the ML equations or prior distribution.  This approach proves to be both a convenient and effective method for estimating the coupling strength and applied current of two reciprocally coupled ML neurons.

\section{Further Research}
These results cover a small fraction of the potential research that can be done with the ML system and parameter approximation. The problem can become increasingly more complex with the alteration of initial conditions, addition of dynamical noise, or estimation of additional parameters.

Extension of this method to higher number of parameters is an obvious goal.  Triangulation of multi-dimensional parameter regions is already possible utilizing N-D Delaunay triangulation via MATLAB's {\sffamily{delaunayn}} and {\sffamily{tsearchn}} functions. The challenge therefore is not in determining membership of a proposed candidate parameter set within an ``NO ROLLBACK'' region, but rather the higher dimensional parameter continuations need to determine the region to be triangularized.  In particular, the current version of XPPAUT does not easily permit continuation beyond 2D, however, the current version of the AUTO library upon which it is based is developing reliable higher dimensional continuation.  Automating the process of LP and Hopf continuation, possibly through scripting, would obviously be advantageous as well. 
 
\begin{center}
\begin{figure}[H]
\centering
\includegraphics[scale=0.25]{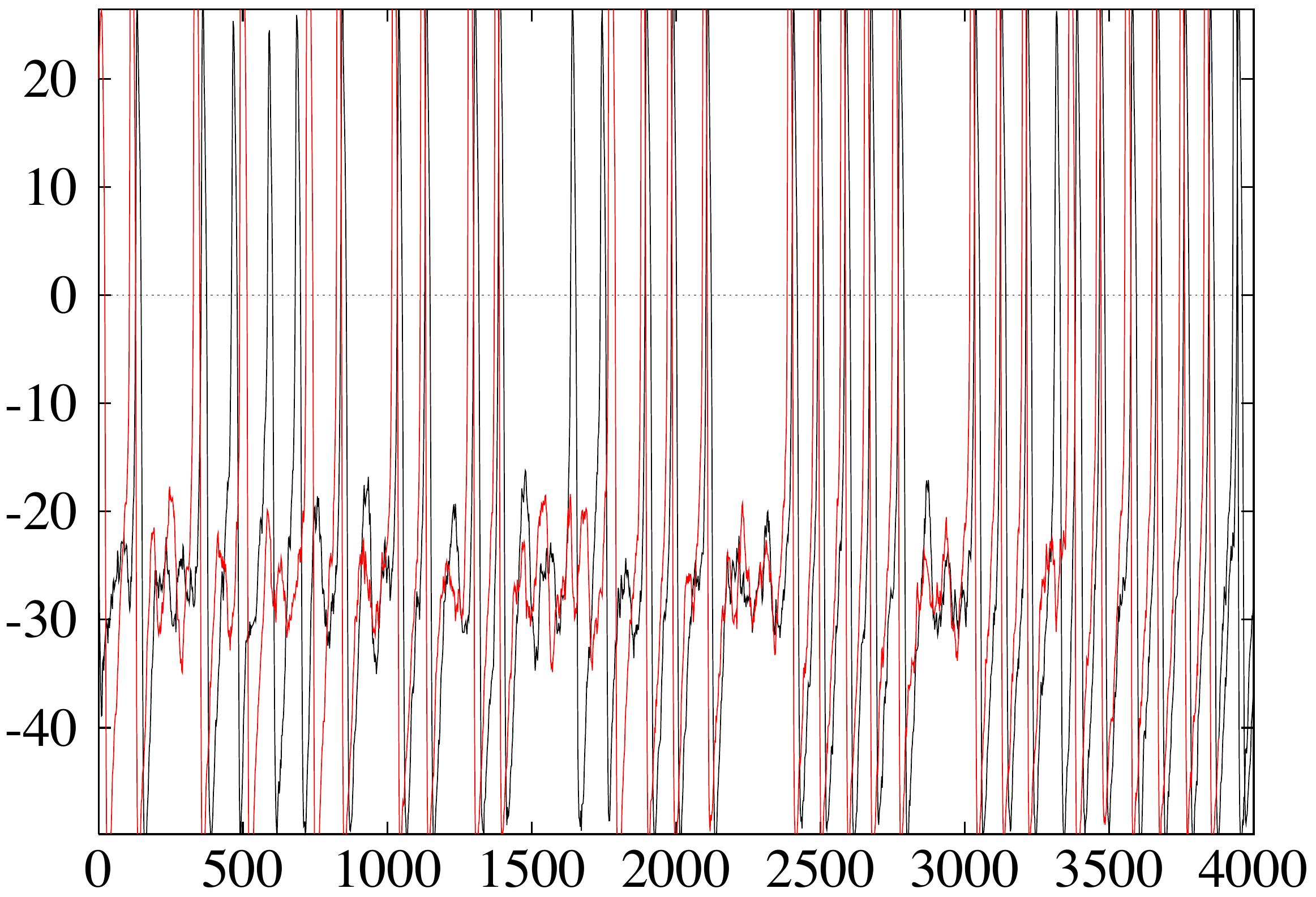}\\
\titlecaption[Mixed mode oscillation]{In the presence of small noise ($\delta_{1}=\delta_{2}=0.7$) and weak coupling ($\gsyn=0.15$) the reciprocally coupled ML model with $\Iapp=95$ predicts low amplitude oscillations disrupted by large amplitude excursions.}
\label{MMO}
\end{figure}
\end{center}

In real-world instances of neuron coupling, interference from nearby neuron firings is present in the form of dynamical noise. This can be accounted for by setting $\delta>0$, in the ML model. The presence of dynamical noise can cause the system to experience MMO, or the stochastic switching between states (shown in Fig.~\ref{MMO}). The parameters of this ML variant may still be able to be estimated, but alterations in the estimation strategy will be necessary. For instance, $\delta$ could be estimated in addition to \Iapp and \gsyn.  Suitable adjustments to the rejection region could be made to account for stochastic perturbations of \Iapp.  

The parameters of the ML system may also be estimated by  conditioning on other statistics. In this study we used cumulative power and mean of the voltage to run MCMC, but this could be altered. For example, conditioning on the difference of cumulative powers of the two voltages was considered, but has not been implemented. This approach may better account for synchrony, or lack thereof, between the two neurons than calculating the likelihoods for each neuron independently.     

\begin{center}
\huge \textbf{Appendix}
\end{center}
\addcontentsline{toc}{section}{Appendix}
\appendix
\refstepcounter{section}

\subsection{Initial Conditions and Parameters}
\label{STDPARAMS}

Unless indicated otherwise, all tests run on the MCMC system had the following initial conditions and parameter values. A biological interpretation of these parameters can be found in Table~\ref{tableofparams}. 
\begin{table}[H]
\centering
\begin{tabular}{|c | c || c | c |}
	\hline
	\textbf{Variable} & \textbf{Initial Value} & \textbf{Variable} & \textbf{Initial Value}  \\ \hline \hline
	$v_1$ & -20 mV & $v_{\text{K}}$& -84 mV \\ \hline
	$v_2$ & 20 mV & $v_{\text{L}}$ & -60 mV\\ \hline
	$w_1$ & 0.3 & $v_{\text{syn}}$& 70 mV \\  \hline
	$w_2$ & 0.5 & $v_{11}$ & -1.2 mV\\  \hline
	$s_1$ & 0.2 & $v_{22}$ & 18 mV \\ \hline
	$s_2$ & 0.1 & $v_3$& 2.0 mV \\ \hline
	C & 20 $\mu$F\ cm$^2$ &$v_4$ & 30 mV \\ \hline
	$g_{\text{Ca}}$ & 4.0 mS/cm$^2$ & $\tau$ & 8 ms \\ \hline
	$g_{\text{K}}$& 8.0 mS/cm$^2$ & $v_{\text{s}}$ & 5.0 mV \\ \hline
	$g_{\text{L}}$ & 2.0 mS/cm$^2$ & $v_{\text{t}}$ & 15 mV\\ \hline
	$v_{\text{Ca}}$ & 120 mS/cm$^2$ & $\phi$ & 0.04/ms \\ \hline
\end{tabular}
\end{table}

\subsection{Cumulative Power of Periodic Functions is \Ot}

Let $m(t)$ have Fourier series representation,\\
	\[m(t)=\frac{a_0}{2}+\sum_{k=1}^n a_k\cos(2\pi \phi_kt)+b_k\sin(2\pi \phi_kt)\]\\
then, it it to be shown that
	\[P(t)=C\cdot{t}+\Oh{1}\]
Proof is By Induction and is adapted from~\citet{quinn}. In the base case ($n=1$),
\begin{align*}\setlength{\arraycolsep}{0pt}
P(t)=\left[8\pi^4 \phi_1^4(a_1^2+b_1^2)\right]t &\overbrace{+2\pi^3 a_1^2 \phi_1^3 \sin(4\pi\phi_1 t)}^{g(t)~\in{~\Oh{1}}}\\
&- 2\pi^3 b_1^2 \phi_1^3 \sin(4\pi\phi_1t)\\
&+16\pi^4a_1b_1\phi_1^4 \cos^2(2\pi\phi_1t)\\
&-16\pi a_1b_1\phi_1^4
\end{align*}
Then it is supposed that for $n-1$ we have,
\begin{align*}\setlength{\arraycolsep}{0pt}
P(t)=\int_0^{t}\left(\sum_{k=1}^{n-1} 4\pi^2 a_k \phi_k^2 \cos(2\pi \phi_kt)+4\pi^2b_k \phi_k^2 \sin(2\pi \phi_kt)\right)^2dt=\left[8\pi^4 \sum_{k=1}^{n-1} \phi_k^4(a_k^2+b_k^2)\right]t+g_{n-1}(t)
\end{align*}
\flushleft{and it remains to be shown that}
\begin{align*}
P(t)=\int_0^{t}\left(\sum_{k=1}^n 4\pi^2 a_k \phi_k^2 \cos(2\pi \phi_kt)+4\pi^2b_k \phi_k^2 \sin(2\pi \phi_kt)\right)^2dt=\left[8\pi^4 \sum_{k=1}^n \phi_k^4(a_k^2+b_k^2)\right]t+g_n(t)
\end{align*}
Next, collecting the $n^{th}$ terms from the summation yields,
\begin{align*}\setlength{\arraycolsep}{0pt}
P(t)=&\int_0^{t}\left(\sum_{k=1}^{n-1} 4\pi^2 a_k \phi_k^2 \cos(2\pi \phi_kt)\right.\\
&+4\pi^2b_k \phi_k^2 \sin(2\pi \phi_kt)+4\pi^4a_n \phi_n^2\cos(2\pi\phi_nt)\\
&\left.+4\pi^4 b_n\phi_n^2\sin(2\pi\phi_nt)\hspace{-3ex}\phantom{\sum_{k=1}^{n-1}}\right)^2 dt\\
\end{align*}
Expanding the square in the previous result gives,
\begin{align*}\setlength{\arraycolsep}{0pt}
&P(t)=\int_0^{t}\left(\sum_{k=1}^{n-1} 4\pi^2 a_k \phi_k^2 \cos(2\pi \phi_kt)+4\pi^2b_k \phi_k^2 \sin(2\pi \phi_kt)\right)^2dt\\
&+2\int_0^t\left(\sum_{k=1}^n 4\pi^2 a_k \phi_k^2 \cos(2\pi \phi_kt)+4\pi^2b_k \phi_k^2 \sin(2\pi \phi_kt)\right)\\
&\quad\quad\cdot\left(4\pi^2a_n\phi_n^2\cos(2\pi\phi_nt)+4\pi^2b_n\phi_n^2\sin(2\pi\phi_nt)\right)dt\\
&+\int_0^t \left(4\pi^2a_n\phi_n^2\cos(2\pi\phi_nt)+4\pi^2b_n\phi_n^2\sin(2\pi\phi_nt)\right)^2dt
\end{align*}
The first term is recognized as the induction hypothesis and so has the form $C\cdot{t}+g_{n-1}(t)$. 
The integrals of the remaining terms are evaluated with the extensive use of trigonometric identities.
\tiny
\begin{align*}\tiny
P(t)=&\left[8\pi^4 \sum_{k=1}^{n-1} \phi_k^4(a_k^2+b_k^2)\right]t+g_{n-1}(t)\\
&+\sum_{k=1}^{n-1}8\pi^3a_ka_n\phi_k^2\phi_n^2\left[\frac{\sin(2\pi(\phi_k-\phi_n)t)}{\phi_k-\phi_n}+\frac{\sin(2\pi(\phi_k+\phi_n)t)}{\phi_k+\phi_n}\right]\\
&-\sum_{k=1}^{n-1}8\pi^3b_ka_n\phi_k^2\phi_n^2\left[\frac{\cos(2\pi(\phi_k-\phi_n)t)}{\phi_k-\phi_n}+\frac{\cos(2\pi(\phi_k+\phi_n)t)}{\phi_k+\phi_n}\right]\\
&+\sum_{k=1}^{n-1}8\pi^3b_ka_n\phi_k^2\phi_n^2\left[\frac{1}{\phi_k-\phi_n}+\frac{1}{\phi_k+\phi_n}\right]\\
&-\sum_{k=1}^{n-1}8\pi^3a_kb_n\phi_k^2\phi_n^2\left[\frac{\cos(2\pi(\phi_k-\phi_n)t)}{\phi_k-\phi_n}+\frac{\cos(2\pi(\phi_k+\phi_n)t)}{\phi_k+\phi_n}\right]\\
&+\sum_{k=1}^{n-1}8\pi^3a_kb_n\phi_k^2\phi_n^2\left[\frac{1}{\phi_k-\phi_n}+\frac{1}{\phi_k+\phi_n}\right]\\
&+\sum_{k=1}^{n-1}8\pi^3b_kb_n\phi_k^2\phi_n^2\left[\frac{\sin(2\pi(\phi_k-\phi_n)t)}{\phi_k-\phi_n}+\frac{\sin(2\pi(\phi_k+\phi_n)t)}{\phi_k+\phi_n}\right]\\
&\mathbf{+8\pi^4\phi_n^4a_n^2t}+2\pi^3a_n^2\phi_n^3\sin(4\pi\phi_nt)\\
&+16\pi^4a_nb_n\phi_n^4\sin^2(2\pi\phi_nt)\\
&\mathbf{+8\pi^4\phi_n^4b_n^2t}-2\pi^3b_n^2\phi_n^3\sin(4\pi\phi_nt)
\end{align*}\normalsize
The bold terms may be combined with the leading term of the induction hypothesis raising the upper bound of the summation from $n-1$ to $n$. The remaining terms, only containing $t$ as arguments of sines and cosines, can be merged with $g_{n-1}(t)$ from the induction hypothesis. Calling this merger $g_n(t)$ completes the induction. Cumulative power has been written in the desired form
\begin{eqnarray*}
P(t)&=&\left[8\pi^4\sum_{k=1}^n \phi_k^4(a_k^2+b_k^2)\right]t+g_n(t)\\
&=&C\cdot{t} + \Oh{1}
\end{eqnarray*}
Since clearly $C\cdot{t}\in\Ot$ \underline{and} $g_n(t)\in\Ot$ it follows that their sum $P(t)\in\Ot$ proving the desired result.

\section*{Acknowledgments} \addcontentsline{toc}{section}{Acknowledgments}
This report summarizes work that was done as part of the Summer Undergraduate Research Institute 
of Experimental Mathematics (SURIEM) held at the Lyman Briggs College of Michigan State University.  
We are very grateful to the National Security Agency and the National Science Foundation for funding 
this research. We would also like to thank our advisor, Professor ~Daniel P. Dougherty, for his guidance 
throughout the summer, and our graduate assistant, Joseph E. ~Roth, for his assistance.


\end{document}